\documentclass[12pt,oneside,reqno]{amsart}
\usepackage{indentfirst}
\usepackage{graphicx}
\usepackage{mathrsfs}
\usepackage{stmaryrd}
\usepackage{amsfonts}
\usepackage{cite}
\usepackage{enumerate,amsmath,amssymb,amsthm}
\usepackage{color}
\usepackage{float}
\newcommand{\dif}{\mathrm{d}}

\newcommand{\be}{\begin{eqnarray}}
\newcommand{\ee}{\end{eqnarray}}
\newcommand{\ce}{\begin{eqnarray*}}
\newcommand{\de}{\end{eqnarray*}}
\newtheorem{theorem}{Theorem}[section]
\newtheorem{lemma}[theorem]{Lemma}
\newtheorem{remark}[theorem]{Remark}
\newtheorem{definition}[theorem]{Definition}
\newtheorem{proposition}[theorem]{Proposition}
\newtheorem{Examples}[theorem]{Examples}
\newtheorem{corollary}[theorem]{Corollary}
\def\e{\varepsilon}
\def\s{\sigma}
\def\t{\theta}

\def\d{\delta}
\def\p{\partial}
\def\g{\gamma}

\def\[{{\Big[}}
\def\]{{\Big]}}
\def\<{{\langle}}
\def\>{{\rangle}}
\def\({{\Big(}}
\def\){{\Big)}}

\def\no{\nonumber}
\def\bt{\begin{theorem}}
\def\et{\end{theorem}}
\def\bl{\begin{lemma}}
\def\el{\end{lemma}}
\def\br{\begin{remark}}
\def\er{\end{remark}}
\def\bx{\begin{Examples}}
\def\ex{\end{Examples}}
\def\bd{\begin{definition}}
\def\ed{\end{definition}}
\def\bp{\begin{proposition}}
\def\ep{\end{proposition}}
\def\bc{\begin{corollary}}
\def\ec{\end{corollary}}

\def\cD{{\mathcal D}}

\def\cK{{\mathcal K}}

\def\cO{{\mathcal O}}

\def\cS{{\mathcal S}}

\def\mE{{\mathbb E}}
\def\mF{{\mathbb F}}

\def\mN{{\mathbb N}}

\def\mP{{\mathbb P}}

\def\mR{{\mathbb R}}
\def\mS{{\mathbb S}}

\def\sA{{\mathscr A}}

\def\sF{{\mathscr F}}

\def\sL{{\mathscr L}}

\def\sV{{\mathscr V}}

\def\geq{\geqslant}
\def\leq{\leqslant}
\allowdisplaybreaks[3]

\begin{document}
	
\allowdisplaybreaks
	
\title{Probabilistic approach to homogenization for a type of multivalued Dirichlet-Neumann problems}
	
\author{Huijie Qiao}
	
	
	
\thanks{This work was supported by NSF of China (No.12071071)}
		
	
\date{}
\dedicatory{School of Mathematics,
		Southeast University\\
		Nanjing, Jiangsu 211189, China\\
		hjqiaogean@seu.edu.cn}

\begin{abstract}
The work is about homogenization for a type of multivalued Dirichlet-Neumann problems. First, we prove an average principle for general multivalued stochastic differential equations in the weak sense. Then for general forward-backward coupled multivalued stochastic systems, the other average principle is presented. Finally, we apply the result to a type of multivalued Dirichlet-Neumann problems and investigate its homogenization.
\end{abstract}

\maketitle \rm

\section{Introduction}
The earliest study about homogenization of partial differential equations (PDEs for short) by probabilistic methods can be traced back to \cite{f}. Then, Bensoussan, Lions and Papanicolaou \cite{blp} elaborated homogenization of a kind of linear parabolic PDEs through stochastic differential equations (SDEs for short). Later, Hu and Peng \cite{hp} studied homogenization of a kind of semilinear parabolic PDEs via forward-backward coupled SDEs. From then on, there have been a lot of results about homogenization for all kinds of semilinear parabolic PDEs, see \cite{bh, bhp, bi, d, do1, do2, eo, l1, l2, p, op} and the references therein. Recently, in \cite{hjw} Hu, Jiang and Wang observed homogenization of a kind of fully nonlinear parabolic PDEs via forward-backward coupled SDEs driven by $G$-Brownian motions.

In this paper, we develop the above mentioned probabilistic approach to homogenization for a type of multivalued Dirichlet-Neumann problems. Concretely speaking, consider the following system with the Dirichlet-Neumann boundary condition
\be\left\{\begin{array}{l}
\frac{\p u^\e(t,x)}{\p t}+(\sL^\e u^\e)(t,x)+f(\frac{t}{\e},x,u^\e(t,x))\in \p \varphi(u^\e(t,x)), ~(t,x)\in(0,T]\times\cO, \\ 
\frac{\p u^\e(t,x)}{\p \Gamma}+g(t,x,u^\e(t,x))\in \p\psi(u^\e(t,x)), ~(t,x)\in[0,T]\times\p \cO,\\
u^\e(T,x)=\Phi(x), \quad x\in\bar{\cO},
\end{array}
\right.
\label{dnpdein}
\ee
where $\varphi, \psi: \mR^d\rightarrow (-\infty, +\infty]$ are convex lower semicontinuous functions, whose subdifferential operators $\p \varphi, \p\psi$ are maximal monotone operators (See Subsection \ref{mmo}), $\sL^\e$ is a second order differential operator defined by
\be
(\sL^\e F)(t,x)=b^i(\frac{t}{\e},x)\p_i F(x)+\frac{1}{2}a^{ij}(\frac{t}{\e},x) \p_i\p_j F(x), \quad F\in C^2(\mR^m),
\label{ledefi}
\ee
and $b: \mR_+\times\mR^{m}\rightarrow\mR^{m}, \s: \mR_+\times\mR^{m}\rightarrow\mR^{m\times m}, f: \mR_+\times\overline{\cO}\times\mR^d\rightarrow\mR^d, g: \mR_+\times\p\cO\times\mR^d\rightarrow\mR^d, \Phi: \overline{\cO}\rightarrow\mR^d$ are all Borel measurable with $a(t,x):=(\s\s^*)(t,x)$. Here $\cO$ is an open connected bounded subset of $\mR^m$ of the form $\cO:=\{x\in\mR^m, \phi(x)>0\}$ with $\p\cO=\{x\in\mR^m, \phi(x)=0\}$, where $\phi\in C_b^3(\mR^m)$ and $|\triangledown \phi(x)|=1$ for $x\in\p\cO$. The internal normal derivative of a function $v\in C^1(\p \cO)$ is given by $\frac{\p v(x)}{\p \Gamma}=-\<\triangledown \phi(x), \triangledown v(x)\>$ for all $x\in\p \cO$. $0<\e<1$ is a small parameter. Under suitable assumptions we observe the convergence of solutions $u^\e(t,x)$ for the system (\ref{dnpdein}) as $\e$ tends to $0$. 

In order to study homogenization for the system (\ref{dnpdein}), we will face two difficulties. One difficulty lies in the PDE itself. Since the operator $\p \varphi$ is nonlinear and non-smooth, we can not straightly describe the limit of $\p \varphi(u^\e(t,x))$. The other difficulty is from the Neumann boundary condition. As $\e\rightarrow 0$, we can not determine whether $\p \psi(u^\e(t,x))$ converges. To overcome these two difficulties, we make use of the link between the system (\ref{dnpdein}) and the following forward-backward coupled multivalued stochastic system, 
\ce\left\{\begin{array}{l}
\dif X_{s}^{\e,t,x}=\triangledown \phi(X_{s}^{\e,t,x})\dif |K^{\e,t,x}|_t^s+b(\frac{s}{\e},X_{s}^{\e,t,x})\dif s+\s(\frac{s}{\e},X_{s}^{\e,t,x})\dif B_{s},\\
K^{\e,t,x}_s=\int_t^s\triangledown \phi(X_{r}^{\e,t,x})\dif |K^{\e,t,x}|_t^r, \quad |K^{\e,t,x}|_t^s=\int_t^s I_{\{X_{r}^{\e,t,x}\in \p\cO\}}\dif |K^{\e,t,x}|_t^r,\\
X_{t}^{\e,t,x}=x\in \overline{\cO},
\end{array}
\right.
\de
and
\ce\left\{\begin{array}{l}
\dif Y_{s}^{\e,t,x}\in\p\varphi(Y_{s}^{\e,t,x})\dif s+\p\psi(Y_{s}^{\e,t,x})\dif |K^{\e,t,x}|_t^s-f(\frac{s}{\e},X_{s}^{\e,t,x},Y_{s}^{\e,t,x})\dif s\\
\qquad\qquad-g(s,X_{s}^{\e,t,x},Y_{s}^{\e,t,x})\dif |K^{\e,t,x}|_t^s+Z_{s}^{\e,t,x}\dif B_{s},\\
Y_{T}^{\e,t,x}=\Phi(X_{T}^{\e,t,x}),
\end{array}
\right.
\de
and prove two average principles for the latter.

In this work, first of all, we prove an average principle for general multivalued stochastic differential equations in the weak sense. Our conditions are weaker than that in some known results (cf. \cite{nn, xl}). Then for forward-backward coupled multivalued stochastic systems, the other average principle is presented. Since $\int_s^Tg(r,X_{r}^{\e,t,x},Y_{r}^{\e,t,x})\dif |K^{\e,t,x}|_t^r$ appears, where $|K^{\e,t,x}|_t^r$ is a continuous increasing process singular with respect to the Lebesgue measure, we work with the $\cS$-topology on the space $D([t,T],\mR^d)$ of right continuous $\mR^d$-valued functions on $[0,T]$ which have left limits. Finally, we apply the result to a type of multivalued Dirichlet-Neumann problems and obtain its homogenization. Moreover, our result can cover \cite[Theorem 4.1]{eo} and \cite[Theorem 6.4]{q1} in some sense.

This paper is arranged as follows. In the next section, we introduce notations, maximal monotone operators and the $\cS$-topology on $D([0,T],\mR^d)$. In Section \ref{main}, main results are formulated. Proofs of two main theorems are placed in Section \ref{xbarxdiffproo} and \ref{averprinproo}, respectively. In Section \ref{app}, we apply our result to a type of multivalued Dirichlet-Neumann problems.

The following convention will be used throughout the paper: $C$, with or without indices, will denote different positive constants whose values may change from one place to another.

\section{Preliminary}\label{fram}

In this section, we introduce notations and concepts, and recall some results used in the sequel.

\subsection{Notations}\label{nota}

In this subsection, we introduce some notations.

For convenience, we shall use $\mid\cdot\mid$ and $\parallel\cdot\parallel$  for norms of vectors and matrices, respectively. Furthermore, let $\langle\cdot$ , $\cdot\rangle$ denote the scalar product in $\mR^m$. Let $U^*$ denote the transpose of a matrix $U$.

Let $C(\mR^m)$ be the collection of continuous functions on $\mR^m$ and $C^2(\mR^m)$ be the space of continuous functions on $\mR^m$ which have continuous partial derivatives of order up to $2$. $C_b^2(\mR^m)$ stands for the subspace of $C^2(\mR^m)$, consisting of functions whose derivatives up to order 2 are bounded. $C_c^2(\mR^m)$ is the collection of all functions in $C^2(\mR^m)$ with compact supports.

\subsection{Maximal monotone operators}\label{mmo}

In this subsection, we introduce maximal monotone operators. 

For a multivalued operator $A: \mR^m\mapsto 2^{\mR^m}$, where $2^{\mR^m}$ stands for all the subsets of $\mR^m$, set
\ce
&&\cD(A):= \left\{x\in \mR^m: A(x) \ne \emptyset\right\},\\
&&Gr(A):= \left\{(x,y)\in \mR^{2m}:x \in \cD(A), ~ y\in A(x)\right\}.
\de
We say that $A$ is monotone if $\langle x_1 - x_2, y_1 - y_2 \rangle \geq 0$ for any $(x_1,y_1), (x_2,y_2) \in Gr(A)$, and $A$ is maximal monotone if 
$$
(x_1,y_1) \in Gr(A) \iff \langle x_1-x_2, y_1 -y_2 \rangle \geq 0, \forall (x_2,y_2) \in Gr(A).
$$

We give two examples to explain maximal monotone operators.

\bx\label{exmmo1}
For a lower semicontinuous convex function $\psi:\mR^m\mapsto(-\infty, +\infty]$, we assume ${\rm Int}(Dom(\psi))\neq \emptyset$, where $Dom(\psi)\equiv\{x\in\mR^m; \psi(x)<\infty\}$ and $\operatorname{Int}(Dom(\psi))$ is the interior of $Dom(\psi)$. Define the subdifferential operator of the function $\psi$:
$$
\partial\psi(x):=\{y\in\mR^m: \<y,z-x\>+\psi(x)\leq \psi(z), \forall z\in\mR^m\}.
$$
Then $\partial\psi$ is a maximal monotone operator.
\ex

\bx\label{exmmo2}
For a closed convex subset $\mathcal{O}$ of $\mathbb{R}^m$, we suppose $\operatorname{Int}(\mathcal{O})\neq\emptyset$. Define the indicator function of $\mathcal{O}$ as follows:
$$
I_{\mathcal{O}}(x):= \begin{cases}0, & \text { if } x \in \mathcal{O}, \\ 
+\infty, & \text { if } x \notin \mathcal{O}.\end{cases}
$$
The subdifferential operator of $I_{\mathcal{O}}$ is given by
$$
\begin{aligned}
\partial I_{\mathcal{O}}(x) & :=\left\{y \in \mathbb{R}^m:\langle y, z-x\rangle \leq 0, \forall z \in \mathcal{O}\right\} \\
& = \begin{cases}\emptyset, & \text { if } x \notin \mathcal{O}, \\
\{0\}, & \text { if } x \in \operatorname{Int}(\mathcal{O}), \\
\Lambda_x, & \text { if } x \in \partial \mathcal{O},\end{cases}
\end{aligned}
$$
where $\Lambda_x$ is the exterior normal cone at $x$. By simple deduction, we know that $\partial I_{\mathcal{O}}$ is a maximal monotone operator.
\ex

Take any $T>0$ and fix it. Let $\sV_{0}$ be the set of all continuous functions $K: [0,T]\mapsto\mR^m$ with finite variations and $K_{0} = 0$. For $K\in\sV_0$ and $s\in [0,T]$, we shall use $|K|_{0}^{s}$ to denote the variation of $K$ on $[0,s]$
and write $|K|_{TV}:=|K|_{0}^{T}$. Set
\ce
&&\sA:=\Big\{(X,K): X\in C([0,T],\overline{\cD(A)}), K \in \sV_0, \\
&&\qquad\qquad\quad~\mbox{and}~\langle X_{t}-x, \dif K_{t}-y\dif t\rangle \geq 0 ~\mbox{for any}~ (x,y)\in Gr(A)\Big\}.
\de
And about $\sA$ we have the following two results (cf.\cite{cepa2, ZXCH}).

\bl\label{equi}
For $X\in C([0,T],\overline{\cD(A)})$ and $K\in \sV_{0}$, the following statements are equivalent:
\begin{enumerate}[(i)]
	\item $(X,K)\in \sA$.
	\item For any $x, y\in C([0,T],\mR^m)$ with $(x_t, y_t)\in Gr(A)$, it holds that 
	$$
	\left\langle X_t-x_t, \dif K_t-y_t\dif t\right\rangle \geq0.
	$$
	\item For any $(X^{'},K^{'})\in \sA$, it holds that 
	$$
	\left\langle X_t-X_t^{'},\dif K_t-\dif K_t^{'}\right\rangle \geq0.
	$$
\end{enumerate}
\el

\bl\label{inteineq}
Assume that ${\rm Int}(\cD(A))\ne\emptyset$. For any $a\in {\rm Int}(\cD(A))$, there exist $M_1 >0$, and $M_{2},M_{3}\geq0$ such that  for any $(X,K)\in \sA$ and $0\leq s<t\leq T$,
$$
\int_s^t{\left< X_r-a, \dif K_r \right>}\geq M_1\left| K \right|_{s}^{t}-M_2\int_s^t{\left| X_r-a\right|}\dif r-M_3\left( t-s \right) .
$$
\el

The following lemma is from \cite[Proposition 3.4]{z}.

\bl\label{finivarilimi}
Suppose that $\left\{X^n; n \geqslant 1\right\}$ and $\left\{K^n; n \geqslant 1\right\}$ are two sequences of continuous $\mathbb{R}^m$ valued processes, defined on possibly different probability spaces $\left(\Omega^n, \sF^n, \mathbb{P}^n\right)$, and converge in distribution to $X$ and $K$, respectively. If for $n \geqslant 1, K^n$ is of bounded variation a.s. and
$$
\lim _{c \rightarrow \infty} \sup _n \mathbb{P}^n\left(\left|K^n\right|_0^T>c\right)=0,
$$
then $K$ is of bounded variation on $[0, T]$ and $\int_0^T\left\langle X^n_t, \dif K^n_t\right\rangle$ converges in distribution to $\int_0^T\langle X_t, \dif K_t\rangle$.
\el

\subsection{The $\cS$-topology on $D([0,T],\mR^d)$}

In this subsection, we introduce the $\cS$-topology on $D([0,T],\mR^d)$.

First of all, we define the $\cS$-topology (\cite[Definition 2.3]{j}). Here let $BV([0,T],\mR^d)$ be the set of all $\mR^d$-valued functions on $[0,T]$  with finite variations. 
Let $C([0,T],\mR^d)$ be the space of all continuous $\mR^d$-valued functions on $[0,T]$ equipped with the uniformly convergence topology.

\bd 
Let $\{x, x_n: n\in\mN\}\subset D([0,T],\mR^d)$. We say that $x_n$ converges to $x$ in the $\cS$-topology if for each $\varepsilon>0$, there exists a sequence $\{v^\e, v_n^{\varepsilon}: n\in\mN\}\subset BV([0,T],\mR^d)$ such that

$(i)$ $\sup\limits _{t \in[0, T]}\left|x(t)-v^{\varepsilon}(t)\right| \leqslant \varepsilon$ and $\sup\limits _{t \in[0, T]}\left|x_n(t)-v_n^{\varepsilon}(t)\right| \leqslant \varepsilon \quad$ for any $n\in\mN$;

$(ii)$ $\lim\limits_{n\rightarrow\infty}v_n^{\varepsilon}(T)=v^{\varepsilon}(T)$ for all $\e>0$;

$(iii)$ for any $f \in C([0,T],\mR^d)$,
$$
\lim _{n \rightarrow \infty} \int_0^T\left\langle f(t), \dif v_n^{\varepsilon}(t)\right\rangle=\int_0^T\left\langle f(t), \dif v^{\varepsilon}(t)\right\rangle .
$$
\ed

The following criterion of tightness under the $\cS$-topology is from \cite[Theorem A.1]{l2}.

\bl\label{stopotigh}
Let $(\Omega, \sF, \mP)$ be a complete probability space and $\{\sF_{t}\}_{t\geq 0}$ be a filtration. Suppose that $\{Z^\e, \e>0\}$ is a family of stochastic processes in $D([0, T],\mathbb{R}^d)$ satisfying
$$
\sup\limits_{\e>0}\left(\mE\sup\limits _{t \in[0, T]}|Z^\e_t|+CV_{[0,T]}(Z^\e)\right)<\infty,
$$
where $CV_{[0,T]}\left(Z^\e\right)$ denotes the conditional variation of $Z^\e$ on $[0,T]$ defined by
$$
CV_{[0,T]}\left(Z^\e\right):=\sup\sum_{i=1}^{n-1} \mathbb{E}\left[\left|\mathbb{E}\left[Z^\e_{t_{i+1}}-Z^\e_{t_i} \mid \sF_{t_i}\right]\right|\right]<\infty
$$
with ``sup" meaning that the supremum is taken over all partitions of the interval $[0,T]$. Then $\left\{Z^\e\right\}$ is tight in the $\cS$-topology. And there exists a subsequence $\left\{\e_k\right\}_{k\in\mN}$ decreasing to $0$, a process $Z$ in $D([0,T],\mR^d)$, and a countable set $I \subset[0, T)$ satisfying that for any $n\in\mN$ and any $\left\{t_1, \cdots, t_n\right\} \subset[0, T] \backslash I,\left(Z^{\e_k}_{t_1}, \cdots, Z^{\e_k}_{t_n}\right)$ converges in distribution to $\left(Z_{t_1}, \cdots, Z_{t_n}\right)$.
\el

\bl\label{martsto}
Let $\left(X^{\varepsilon}, M^{\varepsilon}\right)$ be a multidimensional process in $D([0, T]; \mathbb{R}^p) (p \in \mathbb{N})$ converging to $(X, M)$ in the $\cS$-topology. Let $\left(\sF_t^{X^{\varepsilon}}\right)_{t \geqslant 0}\left(resp. \left(\sF_t^{X}\right)_{t \geqslant 0}\right)$ be the minimal complete admissible filtration for $X^\e$ (resp. $X$). We assume that
$$
\sup\limits _{\varepsilon>0} \mathbb{E}\left[\sup\limits_{0 \leqslant t \leqslant T}\left|M_t^{\varepsilon}\right|^2\right]<C_T, \quad \forall T>0,
$$
$M^{\varepsilon}$ is a $\mathscr{F}^{X^{\varepsilon}}$-martingale, and $M$ is $\mathscr{F}^{X}$-adapted. Then $M$ is a $\mathscr{F}^{X}$-martingale.
\el

\section{Main results}\label{main}

In this section, we formulate the main results.

\subsection{The average principle for multivalued SDEs}

In this subsection, we give an average principle for multivalued SDEs.

Given a complete probability space $(\Omega, \sF, \mP)$ on which a $m$-dimensional standard Brownian motion $B$ is defined. Assume that $\mF:=(\sF_t)_{t\geq 0}$ is the $\mP$-augmentation for the natural filtration of $B$. Fix $T>0$ and consider the following multivalued SDE: $0\leq t\leq s\leq T$
\be\left\{\begin{array}{l}
\dif X_{s}^{\e,t,\xi}\in-A(X_{s}^{\e,t,\xi})\dif s+b(\frac{s}{\e},X_{s}^{\e,t,\xi})\dif s+\s(\frac{s}{\e},X_{s}^{\e,t,\xi})\dif B_{s},\\
X_{t}^{\e,t,\xi}=\xi\in \overline{\cD(A)},
\end{array}
\right.
\label{multfsde}
\ee
where $A$ is a maximal monotone operator with ${\rm Int}(\cD(A))\neq \emptyset$, $b: \mR_+\times\mR^{m}\rightarrow\mR^{m}$ and $\s: \mR_+\times\mR^{m}\rightarrow\mR^{m\times m}$ are all Borel measurable and $\xi$ is a $\sF_t$-measurable random vector. 

We assume:
\begin{enumerate}[$(\mathbf{H}_{A})$]
\item $0\in{\rm Int}(\cD(A))$.
\end{enumerate}
\begin{enumerate}[$(\mathbf{H}^1_{b, \s})$]
\item There exists a constant $L_1>0$ such that for any $s\in\mR_+, x, x_i\in\mR^m, i=1,2,$
\ce
&&|b(s,x_1)-b(s,x_2)|^2+\|\s(s,x_1)-\s(s,x_2)\|^2\leq L_1|x_1-x_2|^2,\\
&&|b(s,x)|^2\leq L_1(1+|x|^2), \quad \|\s(s,x)\|^2\leq L_1.
\de
\end{enumerate}
\begin{enumerate}[$(\mathbf{H}^2_{b, \s})$]
\item 
There exist $\bar{b}: \mR^{m}\rightarrow\mR^{m}$, $\bar{\s}: \mR^{m}\rightarrow\mS_+(\mR^m)$, where $\mS_+(\mR^m)$ is the set of nonnegative definite symmetric $m\times m$ real matrices, such that for any $x\in\mR^m$
\ce
\lim\limits_{\hat{T}\rightarrow\infty}\frac{1}{\hat{T}}\int_0^{\hat{T}}b(s,x)\dif s=\bar{b}(x),\quad \lim\limits_{\hat{T}\rightarrow\infty}\frac{1}{\hat{T}}\int_0^{\hat{T}}(\s\s^*)(s,x)\dif s=(\bar{\s}\bar{\s})(x).
\de
\end{enumerate}
\begin{enumerate}[$(\mathbf{H}^3_{\s})$]
\item
There exists a constant $\iota>0$ such that 
\ce
\<(\s\s^*)(s,x)h, h\>\geq \iota|h|^2, \quad s\in[0,T], x\in\mR^m, h\in\mR^m.
\de
\end{enumerate}

\br\label{condrema}
$(i)$ $(\mathbf{H}_{A})$ is equivalent to that ${\rm Int}(\cD(A))\neq\emptyset$.

$(ii)$ $\bar{b}, \bar{\s}\bar{\s}$ are Lipschitz continuous. Here we only justify $\bar{\s}\bar{\s}$ and by the same way $\bar{b}$ can be verified. In fact, by $(\mathbf{H}^1_{b, \s})$, it holds that for $x_i\in\mR^m, i=1,2,$
\ce
\|(\bar{\s}\bar{\s})(x_1)-(\bar{\s}\bar{\s})(x_2)\|&=&\left\|\lim\limits_{\hat{T}\rightarrow\infty}\frac{1}{\hat{T}}\int_0^{\hat{T}}(\s\s^*)(s,x_1)\dif s-\lim\limits_{\hat{T}\rightarrow\infty}\frac{1}{\hat{T}}\int_0^{\hat{T}}(\s\s^*)(s,x_2)\dif s\right\|\\
&\leq&\lim\limits_{\hat{T}\rightarrow\infty}\frac{1}{\hat{T}}\int_0^{\hat{T}}\|(\s\s^*)(s,x_1)-(\s\s^*)(s,x_2)\|\dif s\\
&\leq&\lim\limits_{\hat{T}\rightarrow\infty}\frac{1}{\hat{T}}\int_0^{\hat{T}}\|\s(s,x_1)-\s(s,x_2)\|\|\s^*(s,x_1)\|\dif s\\
&&+\lim\limits_{\hat{T}\rightarrow\infty}\frac{1}{\hat{T}}\int_0^{\hat{T}}\|\s(s,x_2)\|\|\s^*(s,x_1)-\s^*(s,x_2)\|\dif s\\
&\leq&2L_1|x_1-x_2|.
\de

$(iii)$ By $(ii)$ and $(\mathbf{H}^3_{\s})$, we conclude that $\bar{\s}$ is Lipschitz continuous. Indeed, note that $(\bar{\s}(x_1)-\bar{\s}(x_2))(\bar{\s}(x_1)+\bar{\s}(x_2))=(\bar{\s}\bar{\s})(x_1)-(\bar{\s}\bar{\s})(x_2)$ for $x_i\in\mR^m, i=1,2$. Thus, one can obtain that
\ce
\|\bar{\s}(x_1)-\bar{\s}(x_2)\|\leq\|(\bar{\s}\bar{\s})(x_1)-(\bar{\s}\bar{\s})(x_2)\|\|(\bar{\s}(x_1)+\bar{\s}(x_2))^{-1}\|\leq C|x_1-x_2|.
\de
\er

Under $(\mathbf{H}_{A})$ and $(\mathbf{H}^1_{b, \s})$, by \cite[Theorem 2.8]{rwz}, Eq.(\ref{multfsde}) has a pathwise unique strong solution $(X^{\e,t,\xi},K^{\e,t,\xi})$. Then we construct the following multivalued SDE:
\be\left\{\begin{array}{l}
\dif \bar{X}_{s}^{t,\xi}\in-A(\bar{X}_{s}^{t,\xi})\dif s+\bar{b}(\bar{X}_{s}^{t,\xi})\dif s+\bar{\s}(\bar{X}_{s}^{t,\xi})\dif B_{s},\\
\bar{X}_{t}^{t,\xi}=\xi\in \overline{\cD(A)}.
\end{array}
\right.
\label{averfsde}
\ee
By Remark \ref{condrema} and \cite[Theorem 2.8]{rwz}, $(\mathbf{H}_{A})$, $(\mathbf{H}^1_{b, \s})$, $(\mathbf{H}^2_{b, \s})$ and $(\mathbf{H}^3_{\s})$ assure that Eq.(\ref{averfsde}) has a pathwise unique strong solution $(\bar{X}^{t,\xi},\bar{K}^{t,\xi})$. The following theorem describes the relationship between Eq.(\ref{multfsde}) and Eq.(\ref{averfsde}).

\bt\label{xbarxdiff}
Suppose that $(\mathbf{H}_{A})$, $(\mathbf{H}^1_{b, \s}), (\mathbf{H}^2_{b, \s})$, $(\mathbf{H}^3_{\s})$ hold and $\mE|\xi|^{2}<\infty$. Then it holds that $(X^{\e,t,\xi},K^{\e,t,\xi})$ converges weakly to $(\bar{X}^{t,\xi},\bar{K}^{t,\xi})$ in $C([t,T],\overline{\cD(A)})\times C([t,T],\mR^m)$ as $\e\rightarrow 0$.
\et

The proof of the above theorem is placed in Section \ref{xbarxdiffproo}.

\subsection{The average principle for backward stochastic variation inequalities}

In this subsection, we present an average principle for backward stochastic variation inequalities (BSVIs for short). 

First of all, we take $A=\p I_{\cO}$, where $\cO$ is an open connected bounded subset of $\mR^m$ of the form $\cO:=\{x\in\mR^m, \phi(x)\geq 0\}$ with $\p\cO=\{x\in\mR^m, \phi(x)=0\}$ for a function $\phi\in C_b^3(\mR^m)$ with $|\triangledown \phi(x)|=1$ for $x\in\p\cO$. Note that at any boundary point $x\in\p\cO$, $\triangledown \phi(x)$ is a unit normal vector to the boundary, pointing towards the interior of $\cO$. Then Eq.(\ref{multfsde}) becomes 
\be\left\{\begin{array}{l}
\dif X_{s}^{\e,t,\xi}=\triangledown \phi(X_{s}^{\e,t,\xi})\dif |K^{\e,t,\xi}|_t^s+b(\frac{s}{\e},X_{s}^{\e,t,\xi})\dif s+\s(\frac{s}{\e},X_{s}^{\e,t,\xi})\dif B_{s},\\
K^{\e,t,\xi}_s=\int_t^s\triangledown \phi(X_{r}^{\e,t,\xi})\dif |K^{\e,t,\xi}|_t^r, \quad |K^{\e,t,\xi}|_t^s=\int_t^s I_{\{X_{r}^{\e,t,\xi}\in \p\cO\}}\dif |K^{\e,t,\xi}|_t^r,\\
X_{t}^{\e,t,\xi}=\xi\in \overline{\cO}.
\end{array}
\right.
\label{multfsdef}
\ee

Consider the following BSVI:
\be\left\{\begin{array}{l}
\dif Y_{s}^{\e,t,\xi}\in\p\varphi(Y_{s}^{\e,t,\xi})\dif s+\p\psi(Y_{s}^{\e,t,\xi})\dif |K^{\e,t,\xi}|_t^s-f(\frac{s}{\e},X_{s}^{\e,t,\xi},Y_{s}^{\e,t,\xi})\dif s\\
\qquad\qquad-g(s,X_{s}^{\e,t,\xi},Y_{s}^{\e,t,\xi})\dif |K^{\e,t,\xi}|_t^s+Z_{s}^{\e,t,\xi}\dif B_{s},\\
Y_{T}^{\e,t,\xi}=\Phi(X_{T}^{\e,t,\xi}),
\end{array}
\right.
\label{multfbsde}
\ee
where $\varphi,\psi$ are two lower semicontinuous convex functions with ${\rm Int}(Dom(\varphi))\neq \emptyset, \\{\rm Int}(Dom(\psi))\neq \emptyset$, and $f: \mR_+\times\overline{\cO}\times\mR^d\rightarrow\mR^d, g: \mR_+\times\p\cO\times\mR^d\rightarrow\mR^d$ and $\Phi: \overline{\cO}\rightarrow\mR^d$ are continuous.

We assume:
\begin{enumerate}[$(\mathbf{H}_{\varphi,\psi})$]
\item $(i)$ $0\in {\rm Int}(Dom(\varphi))\cap{\rm Int}(Dom(\psi))$ such that
$$
\varphi(y)\geq\varphi(0)=0, \quad \psi(y)\geq\psi(0)=0, \quad \forall y\in\mR^d.
$$
$(ii)$ 
\ce
\sup\limits_{x\in\overline{\cO}}|\varphi(\Phi(x))|+\sup\limits_{x\in\p\cO}|\psi(\Phi(x))|<\infty.
\de
$(iii)$ There exists a constant $L_2>0$ such that for all $\g>0, s\in[t,T], y\in\mR^d, z\in\mR^{d\times m}$,
\ce
&&\<\triangledown \varphi_{\g}(y),\triangledown \psi_{\g}(y)\>\geq 0,\\
&&\<\triangledown \varphi_{\g}(y),g(s,x,y)\>\leq L_2|\triangledown \psi_{\g}(y)|(1+|g(s,x,y)|), \quad x\in \p\cO, \\
&&\<\triangledown \psi_{\g}(y),f(s,x,y)\>\leq L_2|\triangledown \varphi_{\g}(y)|(1+|f(s,x,y)|), \quad x\in\overline{\cO}, \\
&&-\<\triangledown \varphi_{\g}(y),g(s,x,0)\>\leq L_2|\triangledown \psi_{\g}(y)|(1+|g(s,x,0)|), \quad x\in \p\cO,\\
&&-\<\triangledown \psi_{\g}(y),f(s,x,0)\>\leq L_2|\triangledown \varphi_{\g}(y)|(1+|f(s,x,0)|), \quad x\in\overline{\cO}, 
\de
where $\varphi_{\g}$ is the Moreau-Yosida approximations of $\varphi$ and $\triangledown \varphi_{\g}(y)$ is the derivative of $\varphi_{\g}(y)$ with respect to $y$, and similarly for $\psi_{\g}, \triangledown \psi_{\g}(y)$.
\end{enumerate}
\begin{enumerate}[$(\mathbf{H}^1_{f,g})$]
\item There exist constants $L_3, L_4>0$ such that for any $s\in\mR_+, y, y_1, y_2\in\mR^d$
\ce
&&|f(s,x_1,y_1)-f(s,x_2,y_2)|^2\leq L_3(|x_1-x_2|^2+|y_1-y_2|^2), \quad x_1,x_2\in\overline{\cO}, \\
&&|f(s,x,y)|^2\leq L_3(1+|x|^2+|y|^2), \quad x\in\overline{\cO}, \\
&&|g(s,x_1,y_1)-g(s,x_2,y_2)|^2\leq L_4(|x_1-x_2|^2+|y_1-y_2|^2),\quad x_1,x_2\in\p\cO,\\
&&|g(s,x,y)|^2\leq L_4(1+|x|^2+|y|^2), \quad x\in\p\cO.
\de
\end{enumerate}
\begin{enumerate}[$(\mathbf{H}^2_{f})$]
\item There exists a $\bar{f}: \overline{\cO}\times\mR^d\rightarrow\mR^d$ such that for any $x\in\overline{\cO}, y\in\mR^d$
\ce
\lim\limits_{\hat{T}\rightarrow\infty}\frac{1}{\hat{T}}\int_0^{\hat{T}}f(s,x,y)\dif s=\bar{f}(x,y).
\de
\end{enumerate}

\br\label{barfcond}
By $(\mathbf{H}^1_{f,g})$ and $(\mathbf{H}^2_{f})$, we know that for $x,x_1,x_2\in\overline{\cO}, y, y_1, y_2\in\mR^d$,
\ce
&&|\bar{f}(x_1,y_1)-\bar{f}(x_2,y_2)|^2\leq L_3(|x_1-x_2|^2+|y_1-y_2|^2),\\
&&|\bar{f}(x,y)|^2\leq 2L_3(1+|x|^2+|y|^2).
\de
\er

Under $({\bf H}_{\varphi,\psi})$ and $(\mathbf{H}^1_{f,g})$, by \cite[Proposition 2.5]{pr1}, the system (\ref{multfbsde}) has a unique solution $(Y^{\e,t,\xi}, Z^{\e,t,\xi}, U^{\e,t,\xi}, V^{\e,t,\xi})$. That is, $(Y^{\e,t,\xi}, Z^{\e,t,\xi}, U^{\e,t,\xi}, V^{\e,t,\xi})$ is a $\mR^d\times\mR^{d\times m}\times\mR^d\times\mR^d$-valued progressively measurable stochastic process such that for all $s\in[t,T], \mP.a.s.$
\be
Y^{\e,t,\xi}_s&=&\Phi(X_{T}^{\e,t,\xi})-\int_s^TU^{\e,t,\xi}_r\dif r-\int_s^TV^{\e,t,\xi}_r\dif |K^{\e,t,\xi}|_t^r+\int_s^Tf(\frac{r}{\e},X_{r}^{\e,t,\xi},Y_{r}^{\e,t,\xi})\dif r\no\\
&&+\int_s^Tg(r,X_{r}^{\e,t,\xi},Y_{r}^{\e,t,\xi})\dif |K^{\e,t,\xi}|_t^r-(M^{\e,t,\xi}_T-M^{\e,t,\xi}_s),
\label{yesatiequa}
\ee
where 
$$
M^{\e,t,\xi}_s:=\int_t^s Z_{r}^{\e,t,\xi}\dif B_{r},
$$
and for all $s_1,s_2\in[t,T], s_1\leq s_2$ and all $v\in\mR^m, \mP$-a.s.
\ce
&&\int_{s_1}^{s_2}\<U^{\e,t,\xi}_r, v-Y_{r}^{\e,t,\xi}\>\dif r+\int_{s_1}^{s_2}\varphi(Y_{r}^{\e,t,\xi})\dif r\leq \int_{s_1}^{s_2}\varphi(v)\dif r,\\
&&\int_{s_1}^{s_2}\<V^{\e,t,\xi}_r, v-Y_{r}^{\e,t,\xi}\>\dif |K^{\e,t,\xi}|_t^r+\int_{s_1}^{s_2}\psi(Y_{r}^{\e,t,\xi})\dif |K^{\e,t,\xi}|_t^r\leq \int_{s_1}^{s_2}\psi(v)\dif |K^{\e,t,\xi}|_t^r.
\de
Moreover, we sometimes consider the extension that $X_{s}^{\e,t,\xi}=\xi, |K^{\e,t,\xi}|_t^s=0, Y_{s}^{\e,t,\xi}=Y_{t}^{\e,t,\xi}, Z_{s}^{\e,t,\xi}=0, U^{\e,t,\xi}_s=0, V^{\e,t,\xi}_s=0$ for $0\leq s<t$ in the sequel. 
 
Next, we construct the following BSVI:
\be\left\{\begin{array}{l}
\dif \bar{Y}_{s}^{t,\xi}\in \p\varphi(\bar{Y}_{s}^{t,\xi})\dif s+\p\psi(\bar{Y}_{s}^{t,\xi})\dif |\bar{K}^{t,\xi}|_t^s-\bar{f}(\bar{X}_{s}^{t,\xi},\bar{Y}_{s}^{t,\xi})\dif s\\
\qquad\qquad-g(s,\bar{X}_{s}^{t,\xi},\bar{Y}_{s}^{t,\xi})\dif |\bar{K}^{t,\xi}|_t^s+\bar{Z}_{s}^{t,\xi}\dif B_{s},\\
\bar{Y}_{T}^{t,\xi}=\Phi(\bar{X}_{T}^{t,\xi}).
\end{array}
\right.
\label{averfbsde}
\ee
Based on Remark \ref{barfcond}, Eq.(\ref{averfbsde}) has a unique solution $(\bar{Y}^{t,\xi},\bar{Z}^{t,\xi},\bar{U}^{t,\xi}, \bar{V}^{t,\xi})$.

Now, it is the position to state the main result in this subsection.

\bt\label{averprin}
Assume that $(\mathbf{H}^{1}_{b, \s})$, $(\mathbf{H}^{2}_{b, \s})$, $(\mathbf{H}^3_{\s})$, $({\bf H}_{\varphi,\psi})$, $(\mathbf{H}^1_{f,g})$, $(\mathbf{H}^{2}_{f})$ hold and $\mE|\xi|^{2}<\infty$. Then it holds that $Y^{\e,t,\xi}$ converges weakly to $\bar{Y}^{t,\xi}$ in $D([t,T],\mR^d)$ equipped with the $\cS$-topology as $\e\rightarrow 0$.
\et

Specially, taking $\xi=x\in\bar{\cO}$, we have the following result.

\bc\label{xaverprin}
Assume that $(\mathbf{H}^{1}_{b, \s})$, $(\mathbf{H}^{2}_{b, \s})$, $(\mathbf{H}^3_{\s})$, $({\bf H}_{\varphi,\psi})$, $(\mathbf{H}^1_{f,g})$ and $(\mathbf{H}^{2}_{f})$ hold. Then $Y^{\e,t,x}_t$ converges to $\bar{Y}^{t,x}_t$ as $\e$ tends to $0$.
\ec

The proofs of Theorem \ref{averprin} and Corollary \ref{xaverprin} are placed in Section \ref{averprinproo}.

\section{Proof of Theorem \ref{xbarxdiff}}\label{xbarxdiffproo}

In this section, our aim is to prove Theorem \ref{xbarxdiff}. We start with a key lemma.

\bl\label{xetigh}
Under the assumption of Theorem \ref{xbarxdiff}, it holds that $\{(X^{\e,t,\xi},K^{\e,t,\xi})\}$ is tight in $C([t,T],\overline{\cD(A)})\times C([t,T],\mR^m)$.
\el
\begin{proof}
We only prove this result for $t=0$. And by the same deduction, one can show this result for $t>0$.

{\bf Step 1.} We prove that
\be
&&\sup\limits_{\e}\mE\sup\limits_{s\in[0,T]}|X_s^{\e,0,\xi}|^2\leq C(1+\mE|\xi|^2), \label{xeboun}\\
&&\sup\limits_{\e}\mE|K^{\e,0,\xi}|_0^T\leq C(1+\mE|\xi|^2),\label{keboun}
\ee
where the constant $C>0$ is independent of $\e$.

Note that $(X^{\e,0,\xi},K^{\e,0,\xi})$ satisfies the following equation
$$
X^{\e,0,\xi}_s=\xi-K_s^{\e,0,\xi}+\int_0^s b(\frac{r}{\e},X_{r}^{\e,0,\xi})\dif r+\int_0^s\s(\frac{r}{\e},X_{r}^{\e,0,\xi})\dif B_{r}.
$$
Thus, by the It\^o formula and Lemma \ref{inteineq}, it holds that
\ce
|X^{\e,0,\xi}_s|^2&=&|\xi|^2-\int_0^s\<X^{\e,0,\xi}_r,\dif K_r^{\e,0,\xi}\>+2\int_0^s\<X^{\e,0,\xi}_r, b(\frac{r}{\e},X_{r}^{\e,0,\xi})\>\dif r\\
&&+2\int_0^s\<X^{\e,0,\xi}_r, \s(\frac{r}{\e},X_{r}^{\e,0,\xi})\dif B_{r}\>+\int_0^s\|\s(\frac{r}{\e},X_{r}^{\e,0,\xi})\|^2\dif r\no\\
&\leq&|\xi|^2-2M_1|K^{\e,0,\xi}|_0^s+2M_2\int_0^s|X^{\e,0,\xi}_r|\dif r+2M_3 s+\int_0^s|X^{\e,0,\xi}_r|^2\dif r\no\\
&&+\int_0^s|b(\frac{r}{\e},X_{r}^{\e,0,\xi})|^2\dif r+2\int_0^s\<X^{\e,0,\xi}_r, \s(\frac{r}{\e},X_{r}^{\e,0,\xi})\dif B_{r}\>\no\\
&&+\int_0^s\|\s(\frac{r}{\e},X_{r}^{\e,0,\xi})\|^2\dif r,
\de
and furthermore by the Burkholder-Davis-Gundy inequality
\be
&&\mE\sup\limits_{s\in[0,T]}|X^{\e,0,\xi}_s|^2+2M_1\mE|K^{\e,0,\xi}|_0^T\no\\
&\leq&\mE|\xi|^2+2(M_2+M_3)T+(2M_2+1)\int_0^T\mE\sup\limits_{s\in[0,r]}|X^{\e,0,\xi}_s|^2\dif r\no\\
&&+\int_0^T\mE|b(\frac{r}{\e},X_{r}^{\e,0,\xi})|^2\dif r+\int_0^T\mE\|\s(\frac{r}{\e},X_{r}^{\e,0,\xi})\|^2\dif r\no\\
&&+2C\mE\left(\int_0^T|X^{\e,0,\xi}_r|^2\|\s(\frac{r}{\e},X_{r}^{\e,0,\xi})\|^2\dif r\right)^{1/2}\no\\
&\leq&\mE|\xi|^2+2(M_2+M_3)T+(2M_2+1)\int_0^T\mE\sup\limits_{s\in[0,r]}|X^{\e,0,\xi}_s|^2\dif r\no\\
&&+\int_0^T\mE|b(\frac{r}{\e},X_{r}^{\e,0,\xi})|^2\dif r+\int_0^T\mE\|\s(\frac{r}{\e},X_{r}^{\e,0,\xi})\|^2\dif r\no\\
&&+\frac{1}{2}\mE\sup\limits_{s\in[0,T]}|X^{\e,0,\xi}_s|^2+C\int_0^T\mE\|\s(\frac{r}{\e},X_{r}^{\e,0,\xi})\|^2\dif r.
\label{itoforexpe}
\ee
So, $(\mathbf{H}^{1}_{b, \s})$ and the Gronwall inequality yield (\ref{xeboun}).

Now, combining (\ref{xeboun}) with (\ref{itoforexpe}), we obtain (\ref{keboun}).

{\bf Step 2.} We prove that for any $\eta>0$
\be
&&\lim\limits_{l\downarrow 0}\sup\limits_{\e}\mP\left\{\sup\limits_{0\leq s\leq u\leq s+l\leq T}|X^{\e,0,\xi}_{u}-X^{\e,0,\xi}_{s}|>\eta\right\}=0,\label{pxeuxes}\\
&&\lim\limits_{l\downarrow 0}\sup\limits_{\e}\mP\left\{\sup\limits_{0\leq s\leq u\leq s+l\leq T}|K^{\e,0,\xi}_{u}-K^{\e,0,\xi}_{s}|>\eta\right\}=0. \label{pkeukes}
\ee

First of all, it holds that for $0\leq s\leq u\leq s+l\leq T$
\be
X^{\e,0,\xi}_{u}-X^{\e,0,\xi}_{s}=-K_{u}^{\e,0,\xi}+K_{s}^{\e,0,\xi}+\int_s^{u} b(\frac{r}{\e},X_{r}^{\e,0,\xi})\dif r+\int_s^{u}\s(\frac{r}{\e},X_{r}^{\e,0,\xi})\dif B_{r}.
\label{xekeexpr}
\ee
The It\^o formula, the H\"older inequality and $(\mathbf{H}^1_{b, \s})$ imply that
\be
&&|X^{\e,0,\xi}_{u}-X^{\e,0,\xi}_{s}|^2\no\\
&=&-2\int_s^{u}\<X^{\e,0,\xi}_{r}-X^{\e,0,\xi}_{s}, \dif K_{r}^{\e,0,\xi}\>+2\int_s^{u}\<X^{\e,0,\xi}_{r}-X^{\e,0,\xi}_{s}, b(\frac{r}{\e},X_{r}^{\e,0,\xi})\>\dif r\no\\
&&+2\int_s^{u}\<X^{\e,0,\xi}_{r}-X^{\e,0,\xi}_{s}, \s(\frac{r}{\e},X_{r}^{\e,0,\xi})\dif B_{r}\>+\int_s^{u}\|\s(\frac{r}{\e},X_{r}^{\e,0,\xi})\|^2\dif r\no\\
&\leq&-2\int_s^{u}\<X^{\e,0,\xi}_{r}-X^{\e,0,\xi}_{s}, \dif K_{r}^{\e,0,\xi}\>+\int_s^{u}|X^{\e,0,\xi}_{r}-X^{\e,0,\xi}_{s}|^2\dif r+L_1\int_s^{u}(1+|X_{r}^{\e,0,\xi}|^2)\dif r\no\\
&&+2\left|\int_s^{u}\<X^{\e,0,\xi}_{r}-X^{\e,0,\xi}_{s}, \s(\frac{r}{\e},X_{r}^{\e,0,\xi})\dif B_{r}\>\right|.
\label{xeuxes}
\ee

Next, we compute $-2\int_s^{u}\<X^{\e,0,\xi}_{r}-X^{\e,0,\xi}_{s}, \dif K_{r}^{\e,0,\xi}\>$. Note that $0\in{\rm Int}(\cD(A))$. Thus, there is a $\t_0>0$ such that for any $R>0$ and $\t<\t_0$,
$$
\left\{x \in B(0,R): d(x,(\overline{\cD(A)})^c) \geqslant \t\right\} \neq \emptyset,
$$
where $B(0,R):=\{x\in\mR^m: |x|\leq R\}$, $d(\cdot,\cdot)$ is the Euclidean distance in $\mR^m$ and $(\overline{\cD(A)})^c$ denotes the complement of $\overline{\cD(A)}$. Set
$$
\rho_R(\t):=\sup \left\{|z|: z \in A(x) \text { for all } x \in B(0,R) \text { with } d\left(x,(\overline{\cD(A)})^c\right) \geqslant \t\right\},
$$
and by the local boundedness of $A$ on ${\rm Int}(\cD(A))$, it holds that
$$
\rho_R(\t)<+\infty.
$$
Again put for any $l>0$
$$
\vartheta_R(l):=\inf \left\{\t \in\left(0, \t_0\right): \rho_R(\t) \leqslant l^{-1 / 2}\right\}, 
$$
and we have that
$$
\rho_R\left(l+\vartheta_R(l)\right) \leqslant l^{-1 / 2} \text { and } \quad \lim _{l \downarrow 0} \vartheta_R(l)=0.
$$
Take $l_R>0$ be such that $l_R+\vartheta_R\left(l_R\right)<\t_0$. For $0<l<l_R \wedge 1$, let $X_s^{\e, 0, \xi, l, R}$ be the projection of $X_s^{\e,0,\xi}$ on $\left\{x \in B(0,R): d\left(x,(\overline{\cD(A)})^c\right) \geqslant l+\vartheta_R(l)\right\}$. Thus, for $Z_s^{\e,0, \xi, l, R} \in A(X_s^{\e,0,\xi, l, R}), \sup\limits_{u\in[0,T]}|X_u^{\e,0,\xi}| \leqslant R$ and $0<u-s<l$, it holds that
\ce
&&-2 \int_s^u\left\langle X_r^{\e,0,\xi}-X_s^{\e,0,\xi}, \dif K_r^{\e,0,\xi} \right\rangle\\
 &=&-2 \int_s^u\left\langle X_r^{\e,0,\xi}-X_s^{\e,0,\xi, l, R}, \dif K_r^{\e,0,\xi} \right\rangle-2 \int_s^u\left\langle X_s^{\e,0,\xi, l, R}-X_s^{\e,0,\xi}, \dif K_r^{\e,0,\xi} \right\rangle \\ 
&\leqslant& -2 \int_s^u\left\langle X_r^{\e,0,\xi}-X_s^{\e,0,\xi, l, R}, Z_s^{\e,0,\xi, l, R}\right\rangle \dif r+2\left(l+\vartheta_R(l)\right)\left|K^{\e,0,\xi} \right|_0^T\\
&\leqslant& 4 l^{1 / 2} R+2\left(l+\vartheta_R(l)\right)\left|K^{\e,0,\xi} \right|_0^T,
\de
and furthermore by (\ref{xeuxes})
\ce
&&\sup _{s \leqslant u\leqslant s+l}\left|X_u^{\e,0,\xi}-X_s^{\e,0,\xi}\right|^2 I_{\{\sup\limits_{u\in[0,T]}|X_u^{\e,0,\xi}| \leqslant R\}} \\
&\leqslant&\left(4 l^{1 / 2} R+2\left(l+\vartheta_R(l)\right)\left|K^{\e,0,\xi} \right|_0^T\right)+\sup _{s \leqslant u \leqslant s+l}\int_s^u|X_{r}^{\e,0,\xi}-X_{s}^{\e,0,\xi}|^2\dif r\\
&&+L_1\sup _{s \leqslant u \leqslant s+l}\int_s^u\left(1+\left|X_r^{\e,0,\xi}\right|^2\right) \dif r\\
&&+C \sup _{s \leqslant u \leqslant s+l}\left|\int_s^{u}\<X^{\e,0,\xi}_{r}-X^{\e,0,\xi}_{s}, \s(\frac{r}{\e},X_{r}^{\e,0,\xi})\dif B_{r}\>\right|I_{\{\sup\limits_{u\in[0,T]}|X_u^{\e,0,\xi}| \leqslant R\}}.
\de

From the above deduction and the Burkholder-Davis-Gundy inequality, it follows that
\ce
&&\mE\sup _{s \leqslant u \leqslant s+l}\left|X_u^{\e,0,\xi}-X_s^{\e,0,\xi}\right|^2 I_{\{\sup\limits_{u\in[0,T]}|X_u^{\e,0,\xi}| \leqslant R\}}\\
&\leq&\left(4 l^{1 / 2} R+2\left(l+\vartheta_R(l)\right)\mE|K^{\e,0,\xi}|_0^T\right)+C(1+\mE|\xi|^{2})l\\
&&+C\mE\left(\int_s^{s+l}|X^{\e,0,\xi}_{r}-X^{\e,0,\xi}_{s}|^2\|\s(\frac{r}{\e},X_{r}^{\e,0,\xi})\|^2 \dif r\right)^{1/2}I_{\{\sup\limits_{u\in[0,T]}|X_u^{\e,0,\xi}| \leqslant R\}}\\
&\leq&\left(4 l^{1 / 2} R+2\left(l+\vartheta_R(l)\right)\mE|K^{\e,0,\xi}|_0^T\right)+C(1+\mE|\xi|^{2})l\\
&&+\frac{1}{2}\mE\sup _{s \leqslant u \leqslant s+l}\left|X_u^{\e,0,\xi}-X_s^{\e,0,\xi}\right|^2 I_{\{\sup\limits_{u\in[0,T]}|X_u^{\e,0,\xi}| \leqslant R\}}+C(1+\mE|\xi|^{2})l,
\de
which together with (\ref{keboun}) yields that
\be
&&\sup\limits_{s\in[0,T]}\mE\sup _{s \leqslant u \leqslant s+l}\left|X_u^{\e,0,\xi}-X_s^{\e,0,\xi}\right|^2 I_{\{\sup\limits_{u\in[0,T]}|X_u^{\e,0,\xi}| \leqslant R\}}\no\\
&\leq& \left(8 l^{1 / 2} R+4\left(l+\vartheta_R(l)\right)C(1+\mE|\xi|^{2})\right)+C(1+\mE|\xi|^{2})l.
\label{xeuses}
\ee
Thus, it holds that
\ce
&&\mP\left\{\sup\limits_{0\leq s\leq u\leq s+l\leq T}|X^{\e,0,\xi}_{u}-X^{\e,0,\xi}_{s}|>\eta\right\}\\
&=&\mP\left\{\sup\limits_{0\leq s\leq u\leq s+l\leq T}|X^{\e,0,\xi}_{u}-X^{\e,0,\xi}_{s}|>\eta,\sup\limits_{u\in[0,T]}|X_u^{\e,0,\xi}| \leqslant R\right\}\\
&&+\mP\left\{\sup\limits_{0\leq s\leq u\leq s+l\leq T}|X^{\e,0,\xi}_{u}-X^{\e,0,\xi}_{s}|>\eta,\sup\limits_{u\in[0,T]}|X_u^{\e,0,\xi}| >R\right\}\\
&\leq&\frac{1}{\eta^2}\mE\sup _{s \leqslant u \leqslant s+l}\left|X_u^{\e,0,\xi}-X_s^{\e,0,\xi}\right|^2 I_{\{\sup\limits_{u\in[0,T]}|X_u^{\e,0,\xi}| \leqslant R\}}+\frac{1}{R^2}\mE\sup\limits_{u\in[0,T]}|X_u^{\e,0,\xi}|^2\\
&\leq&\frac{1}{\eta^2}\[\left(8 l^{1 / 2} R+4\left(l+\vartheta_R(l)\right)C(1+\mE|\xi|^{2})\right)+C(1+\mE|\xi|^{2})l\]+\frac{1}{R^2}C(1+\mE|\xi|^{2}),
\de
where we use the Chebyshev inequality. Letting $l\rightarrow 0$ first and then $R\rightarrow\infty$, one can obtain (\ref{pxeuxes}). So, from (\ref{pxeuxes}) and (\ref{xekeexpr}), (\ref{pkeukes}) follows.

Finally, combining (\ref{pxeuxes}), (\ref{pkeukes}) with \cite[Lemma 20.2]{h}, we conclude that $\{(X^{\e,t,\xi},K^{\e,t,\xi})\}$ is tight in $C([t,T],\overline{\cD(A)})\times C([t,T],\mR^m)$.
\end{proof}

\bl\label{xeconver}
Under the assumption of Theorem \ref{xbarxdiff}, it holds that for any $F\in C_c^2(\mR^m)$ and $t\leq s\leq v\leq T$
\ce
&&\lim\limits_{\e\rightarrow0}\mE\Bigg\{\chi_s(X^{\e,t,\xi}, K^{\e,t,\xi})\bigg[F(X^{\e,t,\xi}_v)-F(X^{\e,t,\xi}_s)-\int_s^v(\bar{\sL}F)(X^{\e,t,\xi}_r)\dif r\\
&&\qquad\qquad -\int_s^v\<\p F(X^{\e,t,\xi}_r), \dif K^{\e,t,\xi}_r\>\bigg]\Bigg\}=0,
\de
where $\chi_s$ is a bounded continuous $\sF_s$-measurable functional, $\bar{a}(x):=(\bar{\s}\bar{\s})(x)$ and the operator $\bar{\sL}$ is given by
\be
(\bar{\sL}F)(x):=\bar{b}^i(x)\p_i F(x)+\frac{1}{2}\bar{a}^{ij}(x)\p_i\p_j F(x).
\label{barldefi}
\ee
\el
\begin{proof}
First of all, applying the It\^o formula to $F(X^{\e,t,\xi}_v)$, we obtain that
\ce
&&F(X^{\e,t,\xi}_v)-F(X^{\e,t,\xi}_s)-\int_s^v (\sL^\e F)(\frac{r}{\e}, X^{\e,t,\xi}_r)\dif r-\int_s^v\<\p F(X^{\e,t,\xi}_r), \dif K^{\e,t,\xi}_r\>\\
&=&\int_s^v\<\p F(X^{\e,t,\xi}_r), \dif B_{r}\> ,
\de
where $a(s,x):=(\s\s^*)(s,x)$ and
\ce
(\sL^\e F)(\frac{r}{\e}, X^{\e,t,\xi}_r):=b^i(\frac{r}{\e}, X^{\e,t,\xi}_r)\p_i F(X^{\e,t,\xi}_r)+\frac{1}{2}a^{ij}(\frac{r}{\e}, X^{\e,t,\xi}_r)\p_i\p_j F(X^{\e,t,\xi}_r).
\de 
Note that $\int_t^{\cdot}\<\p F(X^{\e,t,\xi}_r), \dif B_{r}\>$ is a $(\sF_s)_{s\geq t}$-martingale. Thus, it holds that
\ce
&&\mE\Bigg\{\chi_s(X^{\e,t,\xi}, K^{\e,t,\xi})\bigg[F(X^{\e,t,\xi}_v)-F(X^{\e,t,\xi}_s)-\int_s^v(\sL^\e F)(\frac{r}{\e}, X^{\e,t,\xi}_r)\dif r\\
&&\qquad -\int_s^v\<\p F(X^{\e,t,\xi}_r), \dif K^{\e,t,\xi}_r\>\bigg]\Bigg\}=0.
\de
Based on the above equality, one can have that
\ce
&&\mE\Bigg\{\chi_s(X^{\e,t,\xi}, K^{\e,t,\xi})\bigg[F(X^{\e,t,\xi}_v)-F(X^{\e,t,\xi}_s)-\int_s^v(\bar{\sL}F)(X^{\e,t,\xi}_r)\dif r\\
&&\qquad\qquad -\int_s^v\<\p F(X^{\e,t,\xi}_r), \dif K^{\e,t,\xi}_r\>\bigg]\Bigg\}\\
&=&\mE\Bigg\{\chi_s(X^{\e,t,\xi}, K^{\e,t,\xi})\bigg[\int_s^v(\sL^\e F)(\frac{r}{\e}, X^{\e,t,\xi}_r)\dif r-\int_s^v(\bar{\sL}F)(X^{\e,t,\xi}_r)\dif r\bigg]\Bigg\}\\
&=&\mE\Bigg\{\chi_s(X^{\e,t,\xi}, K^{\e,t,\xi})\bigg[\int_s^v\left(b^i(\frac{r}{\e}, X^{\e,t,\xi}_r)\p_i F(X^{\e,t,\xi}_r)-b^i(\frac{r}{\e}, X^{\e,t,\xi}_{r(\d)})\p_i F(X^{\e,t,\xi}_{r(\d)})\right)\dif r\bigg]\Bigg\}\\
&&+\mE\Bigg\{\chi_s(X^{\e,t,\xi}, K^{\e,t,\xi})\bigg[\int_s^v\left(b^i(\frac{r}{\e}, X^{\e,t,\xi}_{r(\d)})\p_i F(X^{\e,t,\xi}_{r(\d)})-\bar{b}^i(X^{\e,t,\xi}_{r(\d)})\p_i F(X^{\e,t,\xi}_{r(\d)})\right)\dif r\bigg]\Bigg\}\\
&&+\mE\Bigg\{\chi_s(X^{\e,t,\xi}, K^{\e,t,\xi})\bigg[\int_s^v\left(\bar{b}^i(X^{\e,t,\xi}_{r(\d)})\p_i F(X^{\e,t,\xi}_{r(\d)})-\bar{b}^i(X^{\e,t,\xi}_r)\p_i F(X^{\e,t,\xi}_r)\right)\dif r\bigg]\Bigg\}\\
&&+\mE\Bigg\{\chi_s(X^{\e,t,\xi}, K^{\e,t,\xi})\bigg[\int_s^v\left(a^{ij}(\frac{r}{\e}, X^{\e,t,\xi}_r)\p_i\p_j F(X^{\e,t,\xi}_r)-a^{ij}(\frac{r}{\e}, X^{\e,t,\xi}_{r(\d)})\p_i\p_j F(X^{\e,t,\xi}_{r(\d)})\right)\dif r\bigg]\Bigg\}\\
&&+\mE\Bigg\{\chi_s(X^{\e,t,\xi}, K^{\e,t,\xi})\bigg[\int_s^v\left(a^{ij}(\frac{r}{\e}, X^{\e,t,\xi}_{r(\d)})\p_i\p_j F(X^{\e,t,\xi}_{r(\d)})-\bar{a}^{ij}(X^{\e,t,\xi}_{r(\d)})\p_i\p_j F(X^{\e,t,\xi}_{r(\d)})\right)\dif r\bigg]\Bigg\}\\
&&+\mE\Bigg\{\chi_s(X^{\e,t,\xi}, K^{\e,t,\xi})\bigg[\int_s^v\left(\bar{a}^{ij}(X^{\e,t,\xi}_{r(\d)})\p_i\p_j F(X^{\e,t,\xi}_{r(\d)})-\bar{a}^{ij}(X^{\e,t,\xi}_r)\p_i\p_j F(X^{\e,t,\xi}_r)\right)\dif r\bigg]\Bigg\}\\
&=:&I_1+I_2+I_3+I_4+I_5+I_6,
\de
where $r(\d):=[\frac{r-s}{\d}]\d+s$ for any $0<\d<v-s$.

Next, for $|I_1+I_3|$, by $(\mathbf{H}^1_{b, \s})$ and Remark \ref{condrema} $(iii)$, it holds that
\ce
|I_1+I_3|&\leq&C\mE\Bigg\{|\chi_s(X^{\e,t,\xi}, K^{\e,t,\xi})|\int_s^v|X^{\e,t,\xi}_r-X^{\e,t,\xi}_{r(\d)}|\dif r\Bigg\}\\
&\leq&C\int_s^v\mE|X^{\e,t,\xi}_r-X^{\e,t,\xi}_{r(\d)}|\dif s\leq C\int_s^v\(\mE|X^{\e,t,\xi}_r-X^{\e,t,\xi}_{r(\d)}|^2\)^{1/2}\dif r\\
&\leq&C\(\sup\limits_{s\in[t,T]}\mE\sup\limits_{s\leq r\leq s+\d}|X^{\e,t,\xi}_r-X^{\e,t,\xi}_s|^2\)^{1/2}\\
&\leq&C\(\sup\limits_{\e}\sup\limits_{s\in[t,T]}\mE\sup\limits_{s\leq r\leq s+\d}|X^{\e,t,\xi}_r-X^{\e,t,\xi}_s|^2\)^{1/2}.
\de
For $I_2$, we notice that
\ce
&&\int_s^v\left(b^i(\frac{r}{\e}, X^{\e,t,\xi}_{r(\d)})-\bar{b}^i(X^{\e,t,\xi}_{r(\d)})\right)\p_i F(X^{\e,t,\xi}_{r(\d)})\dif r\\
&=&\int_0^{v-s}\left(b^i(\frac{\varrho+s}{\e}, X^{\e,t,\xi}_{[\frac{\varrho}{\d}]\d+s})-\bar{b}^i(X^{\e,t,\xi}_{[\frac{\varrho}{\d}]\d+s})\right)\p_i F(X^{\e,t,\xi}_{[\frac{\varrho}{\d}]\d+s})\dif \varrho\\
&=&\int_0^{[\frac{v-s}{\d}]\d}\left(b^i(\frac{\varrho+s}{\e}, X^{\e,t,\xi}_{[\frac{\varrho}{\d}]\d+s})-\bar{b}^i(X^{\e,t,\xi}_{[\frac{\varrho}{\d}]\d+s})\right)\p_i F(X^{\e,t,\xi}_{[\frac{\varrho}{\d}]\d+s})\dif \varrho\\
&&+\int_{[\frac{v-s}{\d}]\d}^{v-s}\left(b^i(\frac{\varrho+s}{\e}, X^{\e,t,\xi}_{[\frac{\varrho}{\d}]\d+s})-\bar{b}^i(X^{\e,t,\xi}_{[\frac{\varrho}{\d}]\d+s})\right)\p_i F(X^{\e,t,\xi}_{[\frac{\varrho}{\d}]\d+s})\dif \varrho\\
&=&\sum\limits_{k=0}^{[\frac{v-s}{\d}]-1}\int_{k\d}^{(k+1)\d}\left(b^i(\frac{\varrho+s}{\e}, X^{\e,t,\xi}_{k\d+s})-\bar{b}^i(X^{\e,t,\xi}_{k\d+s})\right)\p_i F(X^{\e,t,\xi}_{k\d+s})\dif \varrho\\
&&+\int_{[\frac{v-s}{\d}]\d}^{v-s}\left(b^i(\frac{\varrho+s}{\e}, X^{\e,t,\xi}_{[\frac{\varrho}{\d}]\d+s})-\bar{b}^i(X^{\e,t,\xi}_{[\frac{\varrho}{\d}]\d+s})\right)\p_i F(X^{\e,t,\xi}_{[\frac{\varrho}{\d}]\d+s})\dif \varrho\\
&=&\sum\limits_{k=0}^{[\frac{v-s}{\d}]-1}\e\int_{\frac{k\d+s}{\e}}^{\frac{(k+1)\d+s}{\e}}\left(b^i(\tau, X^{\e,t,\xi}_{k\d+s})-\bar{b}^i(X^{\e,t,\xi}_{k\d+s})\right)\dif \tau\p_i F(X^{\e,t,\xi}_{k\d+s})\\
&&+\int_{[\frac{v-s}{\d}]\d}^{v-s}\left(b^i(\frac{\varrho+s}{\e}, X^{\e,t,\xi}_{[\frac{\varrho}{\d}]\d+s})-\bar{b}^i(X^{\e,t,\xi}_{[\frac{\varrho}{\d}]\d+s})\right)\p_i F(X^{\e,t,\xi}_{[\frac{\varrho}{\d}]\d+s})\dif \varrho.
\de
Thus, from the H\"older inequality and $(\mathbf{H}^1_{b, \s})$, it follows that
\ce
|I_2|&\leq& C\sum\limits_{i=1}^m\mE\sum\limits_{k=0}^{[\frac{v-s}{\d}]-1}\left|\e\int_{\frac{k\d+s}{\e}}^{\frac{(k+1)\d+s}{\e}}\left(b^i(\tau, X^{\e,t,\xi}_{k\d+s})-\bar{b}^i(X^{\e,t,\xi}_{k\d+s})\right)\dif \tau\right|\\
&&+C\mE\left|\int_{[\frac{v-s}{\d}]\d}^{v-s}\left(b^i(\frac{\varrho+s}{\e}, X^{\e,t,\xi}_{[\frac{\varrho}{\d}]\d+s})-\bar{b}^i(X^{\e,t,\xi}_{[\frac{\varrho}{\d}]\d+s})\right)\p_i F(X^{\e,t,\xi}_{[\frac{\varrho}{\d}]\d+s})\dif \varrho\right|\\
&\leq& C\sum\limits_{i=1}^m\mE\sum\limits_{k=0}^{[\frac{v-s}{\d}]-1}\left|\e\int_{\frac{k\d+s}{\e}}^{\frac{(k+1)\d+s}{\e}}\left(b^i(\tau, X^{\e,t,\xi}_{k\d+s})-\bar{b}^i(X^{\e,t,\xi}_{k\d+s})\right)\dif \tau\right|\\
&&+C\sum\limits_{i=1}^m\mE\int_{[\frac{v-s}{\d}]\d}^{v-s}\left|b^i(\frac{\varrho+s}{\e}, X^{\e,t,\xi}_{[\frac{\varrho}{\d}]\d+s})-\bar{b}^i(X^{\e,t,\xi}_{[\frac{\varrho}{\d}]\d+s})\right|\dif \varrho\\
&\leq& C\sum\limits_{i=1}^m\mE\sum\limits_{k=0}^{[\frac{v-s}{\d}]-1}\Bigg|\e\int_{0}^{\frac{(k+1)\d+s}{\e}}\left(b^i(\tau, X^{\e,t,\xi}_{k\d+s})-\bar{b}^i(X^{\e,t,\xi}_{k\d+s})\right)\dif \tau\\
&&\qquad\qquad\qquad\qquad -\e\int_{0}^{\frac{k\d+s}{\e}}\left(b^i(\tau, X^{\e,t,\xi}_{k\d+s})-\bar{b}^i(X^{\e,t,\xi}_{k\d+s})\right)\dif \tau\Bigg|\\
&&+C\d\left(\mE(1+\sup\limits_{r\in[t,T]}|X^{\e,t,\xi}_r|^2)\right)^{1/2}\\
&\leq& C\sum\limits_{i=1}^m\mE\sum\limits_{k=0}^{[\frac{v-s}{\d}]-1}\Bigg|\e\int_{0}^{\frac{(k+1)\d+s}{\e}}\left(b^i(\tau, X^{\e,t,\xi}_{k\d+s})-\bar{b}^i(X^{\e,t,\xi}_{k\d+s})\right)\dif \tau\\
&&\qquad\qquad\qquad\qquad -\e\int_{0}^{\frac{k\d+s}{\e}}\left(b^i(\tau, X^{\e,t,\xi}_{k\d+s})-\bar{b}^i(X^{\e,t,\xi}_{k\d+s})\right)\dif \tau\Bigg|\\
&&+C\d\left(1+\mE|\xi|^{2}\right)^{1/2}.
\de

By the similar deduction to that for $I_1, I_2, I_3$, one can obtain that
\ce
|I_4+I_6|+|I_5|&\leq& C\(\sup\limits_{\e}\sup\limits_{s\in[t,T]}\mE\sup\limits_{s\leq r\leq s+\d}|X^{\e,t,\xi}_r-X^{\e,t,\xi}_s|^2\)^{1/2}\\
&&+C\sum\limits_{i,j=1}^m\mE\sum\limits_{k=0}^{[\frac{v-s}{\d}]-1}\Bigg|\e\int_{0}^{\frac{(k+1)\d+s}{\e}}\left(a^{ij}(\tau, X^{\e,t,\xi}_{k\d+s})-\bar{a}^{ij}(X^{\e,t,\xi}_{k\d+s})\right)\dif \tau\\
&&\qquad\qquad\qquad\qquad -\e\int_{0}^{\frac{k\d+s}{\e}}\left(a^{ij}(\tau, X^{\e,t,\xi}_{k\d+s})-\bar{a}^{ij}(X^{\e,t,\xi}_{k\d+s})\right)\dif \tau\Bigg|\\
&&+C\d\left(1+\mE|\xi|^{2}\right)^{1/2}.
\de

Now, combining the above deduction, by $(\mathbf{H}^2_{b, \s})$ and the dominated convergence theorem we obtain that
\ce
0&\leq& \lim\limits_{\e\rightarrow0}|I_1+I_2+I_3+I_4+I_5+I_6|\\
&\leq& C\(\sup\limits_{\e}\sup\limits_{s\in[t,T]}\mE\sup\limits_{s\leq r\leq s+\d}|X^{\e,t,\xi}_r-X^{\e,t,\xi}_s|^2\)^{1/2}+C\d\left(1+\mE|\xi|^{2}\right)^{1/2}.
\de
Finally, as $\d$ tends to $0$, (\ref{xeuses}) gives the required result.
\end{proof}

At present, we are ready to prove Theorem \ref{xbarxdiff}.

{\bf Proof of Theorem \ref{xbarxdiff}.} By Lemma \ref{xetigh}, there exist a subsequence $\{(X^{\e_k,t,\xi},K^{\e_k,t,\xi})\}$ and a pair process $(\hat{X}^{t,\hat\xi}, \hat{K}^{t,\hat\xi})$ defined on the other complete filtered probability space $(\hat\Omega, \hat\sF, \hat\mP, (\hat\sF_s)_{s\in[t,T]})$ such that as $k$ tends to $\infty$,
\ce
(X^{\e_k,t,\xi},K^{\e_k,t,\xi}) \rightarrow (\hat{X}^{t,\hat\xi}, \hat{K}^{t,\hat\xi}), \quad \mbox{weakly in}~ C([t,T],\overline{\cD(A)})\times C([t,T],\mR^m),
\de
which together with Lemma \ref{finivarilimi} and Lemma \ref{xeconver} yields that
\ce
&&\hat\mE\Bigg\{\chi_s(\hat X^{t,\xi}, \hat K^{t,\xi})\bigg[F(\hat{X}^{t,\xi}_v)-F(\hat{X}^{t,\xi}_s)-\int_s^v(\bar{\sL}F)(\hat{X}^{t,\xi}_r)\dif r\\
&&\qquad\qquad -\int_s^v\<\p F(\hat{X}^{t,\xi}_r), \dif \hat{K}^{t,\xi}_r\>\bigg]\Bigg\}=0,
\de
where $\hat\mE$ denotes the expectation with respect to $\hat\mP$. From this, it follows that $\hat\mP\circ(\hat X^{t,\hat\xi}, \hat K^{t,\hat\xi})^{-1}$ is a solution of the martingale problem for $(\bar{\sL},\mP\circ\xi^{-1})$. Note that the martingale problem for $(\bar{\sL},\mP\circ\xi^{-1})$ is well-posed. Therefore, $\hat\mP\circ(\hat X^{t,\hat\xi}, \hat K^{t,\hat\xi})^{-1}=\mP\circ(\bar X^{t,\xi}, \bar K^{t,\xi})^{-1}$ and $(X^{\e_k,t,\xi},K^{\e_k,t,\xi})$ converges weakly to $(\bar X^{t,\xi}, \bar K^{t,\xi})$ as $k\rightarrow\infty$.

\section{Proofs of Theorem \ref{averprin} and Corollary \ref{xaverprin}}\label{averprinproo}

In this section, we prove Theorem \ref{averprin} and Corollary \ref{xaverprin}.

 First of all, we recall that $(Y^{\e,t,\xi}, Z^{\e,t,\xi}, U^{\e,t,\xi}, V^{\e,t,\xi})$ solves Eq(\ref{yesatiequa}), i.e.
\ce
Y^{\e,t,\xi}_s&=&\Phi(X_{T}^{\e,t,\xi})-\int_s^TU^{\e,t,\xi}_r\dif r-\int_s^TV^{\e,t,\xi}_r\dif |K^{\e,t,\xi}|_t^r+\int_s^Tf(\frac{r}{\e},X_{r}^{\e,t,\xi},Y_{r}^{\e,t,\xi})\dif r\no\\
&&+\int_s^Tg(r,X_{r}^{\e,t,\xi},Y_{r}^{\e,t,\xi})\dif |K^{\e,t,\xi}|_t^r-(M_{T}^{\e,t,\xi}-M_{s}^{\e,t,\xi}), \quad s\in[t,T], \quad \mP.a.s..
\de
Then we investigate the convergence of every part in the above equation.

For any $t\leq s\leq T$, set 
\ce
{\bf U}_s^\e:=\int_t^s U^{\e,t,\xi}_r\dif r, \quad {\bf V}_s^\e:=\int_t^s V^{\e,t,\xi}_r\dif |K^{\e,t,\xi}|_t^r,
\de
and we have the following result.

\bl\label{8proctigh}
Under the assumptions of Theorem \ref{averprin}, $\{(X^{\e,t,\xi}, K^{\e,t,\xi}, Y^{\e,t,\xi}, M^{\e,t,\xi}, {\bf U}^\e, {\bf V}^\e)\}$ is tight in $C([t,T],\overline{\cO})\times C([t,T],\mR^m)\times D([t,T],\mR^d) \times D([t,T],\mR^d) \times D([t,T],\mR^d)\times D([t,T],\mR^d)$  under the $\cS$-topology.
\el
\begin{proof}
First of all, by Lemma \ref{xetigh} and the fact that the $\cS$-topology is weaker than the uniformly convergence topology, we know that $\{(X^{\e,t,\xi}, K^{\e,t,\xi})\}$ is tight under the $\cS$-topology. 

For $Y^{\e,t,\xi}$, by the similar deduction to that for \cite[Theorem 9]{mr1}, it holds that 
\be
&&\sup\limits_{\e}\left[\mE\sup\limits_{r\in[t,T]}|Y^{\e,t,\xi}_r|^2+\mE\int_t^T|Y^{\e,t,\xi}_r|^2(\dif r+\dif |K^{\e,t,\xi}|_t^r)\right]\leq C,\label{yees}\\
&&\sup\limits_{\e}\mE\int_t^T|Z^{\e,t,\xi}_r|^2\dif r\leq C,\label{zees}\\
&&\sup\limits_{\e}\left[\mE\int_t^T|U^{\e,t,\xi}_r|^2\dif r+\mE\int_t^T|V^{\e,t,\xi}_r|^2\dif |K^{\e,t,\xi}|_t^r\right]\leq C, \label{uevees}
\ee
where $C>0$ is a constant depending on $L_3, L_4$. Besides, for any partition $\pi$ of $[t,T]$: $t=t_1<t_2<\cdots<t_n=T$, the H\"older inequality and $(\mathbf{H}^1_{f,g})$ imply that
\ce
&&\sum\limits_{i=1}^{n-1}\mE\left[\left|\mE[Y^{\e,t,\xi}_{t_{i+1}}-Y^{\e,t,\xi}_{t_{i}}|\sF_{t_{i}}]\right|\right]\\
&\leq&\mE\int_t^T|U^{\e,t,\xi}_r|\dif r+\mE\int_t^T|V^{\e,t,\xi}_r|\dif |K^{\e,t,\xi}|_t^r+\mE\int_t^T|f(\frac{r}{\e},X_{r}^{\e,t,\xi},Y_{r}^{\e,t,\xi})|\dif r\\
&&+\mE\int_s^T|g(r,X_{r}^{\e,t,\xi},Y_{r}^{\e,t,\xi})|\dif |K^{\e,t,\xi}|_t^r\\
&\leq&T^{1/2}\left(\mE\int_t^T|U^{\e,t,\xi}_r|^2\dif r\right)^{1/2}+\left(\mE\int_t^T|V^{\e,t,\xi}_r|^2\dif |K^{\e,t,\xi}|_t^r\right)^{1/2}\left(\mE|K^{\e,t,\xi}|_t^T\right)^{1/2}\\
&&+C\mE\int_t^T(1+|X_{r}^{\e,t,\xi}|+|Y_{r}^{\e,t,\xi}|)(\dif r+\dif |K^{\e,t,\xi}|_t^r)\\
&\leq&T^{1/2}\left(\mE\int_t^T|U^{\e,t,\xi}_r|^2\dif r\right)^{1/2}+\left(\mE\int_t^T|V^{\e,t,\xi}_r|^2\dif |K^{\e,t,\xi}|_t^r\right)^{1/2}\left(\mE|K^{\e,t,\xi}|_t^T\right)^{1/2}\\
&&+C\left(\mE\int_t^T(1+|Y_{r}^{\e,t,\xi}|)^2(\dif r+\dif |K^{\e,t,\xi}|_t^r)\right)^{1/2}\left(T+\mE|K^{\e,t,\xi}|_t^T\right)^{1/2}\\
&&+C\left(\mE\sup\limits_{r\in[t,T]}|X^{\e,t,\xi}_r|^2\right)^{1/2}\left(T^2+\mE(|K^{\e,t,\xi}|_t^T)^2\right)^{1/2},
\de
which together with (\ref{uevees}), (\ref{yees}), (\ref{keboun}), (\ref{xeboun}) and \cite[Proposition 1]{mr2} yields that 
$$
\sup\limits_{\e>0}CV_{[t,T]}(Y^{\e,t,\xi})<\infty. 
$$

Now, we deal with $M^{\e,t,\xi}$. By (\ref{zees}), $M^{\e,t,\xi}$ is a martingale with respect to $(\sF_s)_{s\geq t}$, which yields that
$$
\sum\limits_{i=1}^{n-1}\mE\left[\left|\mE[M^{\e,t,\xi}_{t_{i+1}}-M^{\e,t,\xi}_{t_{i}}|\sF_{t_{i}}]\right|\right]=0.
$$ 
Moreover, again by (\ref{zees}) it holds that
\ce
\sup\limits_{\e}\mE\sup\limits_{r\in[t,T]}|M^{\e,t,\xi}_r|\leq C\sup\limits_{\e}\mE\left(\int_t^T |Z_{r}^{\e,t,\xi}|^2\dif r\right)^{1/2}\leq C\left(\sup\limits_{\e}\mE\int_t^T |Z_{r}^{\e,t,\xi}|^2\dif r\right)^{1/2}<\infty.
\de

For ${\bf U}^\e$, by (\ref{uevees}) we have that
\ce
\sup\limits_{\e}\mE\sup\limits_{r\in[t,T]}|{\bf U}^\e_r|\leq \sup\limits_{\e}\mE\int_t^T |U^{\e,t,\xi}_r|\dif r\leq T^{1/2}\left(\sup\limits_{\e}\mE\int_t^T|U^{\e,t,\xi}_r|^2\dif r\right)^{1/2}<\infty,
\de
and 
\ce
\sup\limits_{\e}CV_{[t,T]}({\bf U}^\e)&=&\sup\limits_{\e}\sup\limits_{\pi}\sum\limits_{i=1}^{n-1}\mE\left[\left|\mE[{\bf U}^\e_{t_{i+1}}-{\bf U}^\e_{t_{i}}|\sF_{t_{i}}]\right|\right]\\
&\leq& \sup\limits_{\e}\mE\int_t^T |U^{\e,t,\xi}_r|\dif r\leq T^{1/2}\left(\sup\limits_{\e}\mE\int_t^T|U^{\e,t,\xi}_r|^2\dif r\right)^{1/2}<\infty.
\de
Moreover, by the similar deduction to that for ${\bf U}^\e$, one can obtain that 
$$
\sup\limits_{\e}\left[\mE\sup\limits_{r\in[t,T]}|{\bf V}^\e_r|+CV_{[t,T]}({\bf V}^\e)\right]<\infty.
$$

Finally, collecting the above estimates, by Lemma \ref{stopotigh} we conclude that $\{(Y^{\e,t,\xi}, M^{\e,t,\xi},\\ {\bf U}^\e, {\bf V}^\e)\}$ is tight with respect to the $\cS$-topology. 
The proof is complete.
\end{proof}

To study the convergence of $\int_s^Tf(\frac{r}{\e},X_{r}^{\e,t,\xi},Y_{r}^{\e,t,\xi})\dif r$, we prepare the following limit result.

\bl\label{yessles}
Under the assumptions of Theorem \ref{averprin}, it holds that
\ce
\lim\limits_{l\rightarrow0}\sup\limits_{\e}\sup\limits_{s\in[t,T]}\sup _{s \leqslant r \leqslant s+l}\mE|Y_{r}^{\e,t,\xi}-Y_{s+l}^{\e,t,\xi}|^{2}=0.
\de
\el
\begin{proof}
First of all, we have that for $l>0$ and $t\leq s\leq r\leq s+l\leq T$
\ce
Y^{\e,t,\xi}_r-Y^{\e,t,\xi}_{s+l}&=&-\int_r^{s+l}U^{\e,t,\xi}_v\dif v-\int_r^{s+l}V^{\e,t,\xi}_v\dif |K^{\e,t,\xi}|_t^v+\int_r^{s+l}f(\frac{v}{\e},X_{v}^{\e,t,\xi},Y_{v}^{\e,t,\xi})\dif v\no\\
&&+\int_r^{s+l}g(v,X_{v}^{\e,t,\xi},Y_{v}^{\e,t,\xi})\dif |K^{\e,t,\xi}|_t^v-\int_r^{s+l}Z_{v}^{\e,t,\xi}\dif B_{v}.
\de
Then the It\^o formula implies that
\be
&&|Y^{\e,t,\xi}_r-Y^{\e,t,\xi}_{s+l}|^2+\int_r^{s+l}\|Z_{v}^{\e,t,\xi}\|^2\dif v\no\\
&=&-2\int_r^{s+l}\<Y^{\e,t,\xi}_v-Y^{\e,t,\xi}_{s+l}, U^{\e,t,\xi}_v\>\dif v-2\int_r^{s+l}\<Y^{\e,t,\xi}_v-Y^{\e,t,\xi}_{s+l},V^{\e,t,\xi}_v\>\dif |K^{\e,t,\xi}|_t^v\no\\
&&+2\int_r^{s+l}\<Y^{\e,t,\xi}_v-Y^{\e,t,\xi}_{s+l}, f(\frac{v}{\e},X_{v}^{\e,t,\xi},Y_{v}^{\e,t,\xi})\>\dif v\no\\
&&+2\int_r^{s+l}\<Y^{\e,t,\xi}_v-Y^{\e,t,\xi}_{s+l}, g(v,X_{v}^{\e,t,\xi},Y_{v}^{\e,t,\xi})\>\dif |K^{\e,t,\xi}|_t^v\no\\
&&-2\int_r^{s+l}\<Y^{\e,t,\xi}_v-Y^{\e,t,\xi}_{s+l}, Z_{v}^{\e,t,\xi}\dif B_{v}\>.
\label{itoye}
\ee

From the above equality and $(\mathbf{H}^1_{f,g})$, it follows that for any $R>0$
\ce
&&|Y^{\e,t,\xi}_r-Y^{\e,t,\xi}_{s+l}|^2I_{\{\sup\limits_{s\in[t,T]}|Y^{\e,t,\xi}_s|<R\}}+\int_r^{s+l}\|Z_{v}^{\e,t,\xi}\|^2\dif vI_{\{\sup\limits_{s\in[t,T]}|Y^{\e,t,\xi}_s|<R\}}\\
&\leq& 4R\int_r^{s+l}|U^{\e,t,\xi}_v|\dif v+4R\int_r^{s+l}|V^{\e,t,\xi}_v|\dif |K^{\e,t,\xi}|_t^v\\
&&+4RL_3^{1/2}\int_r^{s+l}(1+|X_{v}^{\e,t,\xi}|+|Y_{v}^{\e,t,\xi}|)\dif v\\
&&+4RL_4^{1/2}\int_r^{s+l}(1+|X_{v}^{\e,t,\xi}|+|Y_{v}^{\e,t,\xi}|)\dif |K^{\e,t,\xi}|_t^v\\
&&-2\int_r^{s+l}\<Y^{\e,t,\xi}_v-Y^{\e,t,\xi}_{s+l}, Z_{v}^{\e,t,\xi}\dif B_{v}\>I_{\{\sup\limits_{s\in[t,T]}|Y^{\e,t,\xi}_s|<R\}}.
\de
Taking the expectation on two sides, by the H\"older inequality we have that
\ce
&&\mE|Y^{\e,t,\xi}_r-Y^{\e,t,\xi}_{s+l}|^2I_{\{\sup\limits_{s\in[t,T]}|Y^{\e,t,\xi}_s|<R\}}+\mE\int_r^{s+l}\|Z_{v}^{\e,t,\xi}\|^2\dif vI_{\{\sup\limits_{s\in[t,T]}|Y^{\e,t,\xi}_s|<R\}}\\
&\leq& 4R\left(\mE\int_r^{s+l}|U^{\e,t,\xi}_v|^2\dif v\right)^{1/2}(s+l-r)^{1/2}\\
&&+4R\left(\mE\int_r^{s+l}|V^{\e,t,\xi}_v|^2\dif |K^{\e,t,\xi}|_t^v\right)^{1/2}\left(\mE(|K^{\e,t,\xi}|_t^{s+l}-|K^{\e,t,\xi}|_t^r)\right)^{1/2}\\
&&+4RL_3^{1/2}\left(\mE\int_r^{s+l}(1+|X_{v}^{\e,t,\xi}|+|Y_{v}^{\e,t,\xi}|)^2\dif v\right)^{1/2}(s+l-r)^{1/2}\\
&&+4RL_4^{1/2}\left(\mE(1+\sup\limits_{v\in[t,T]}|X_{v}^{\e,t,\xi}|+\sup\limits_{v\in[t,T]}|Y_{v}^{\e,t,\xi}|)^2\right)^{1/2}\(\mE(|K^{\e,t,\xi}|_t^{s+l}-|K^{\e,t,\xi}|_t^r)^2\)^{1/2},
\de
which together with (\ref{uevees}), (\ref{yees}), (\ref{xeboun}) yields that
\ce
&&\sup _{s \leqslant r \leqslant s+l}\mE|Y^{\e,t,\xi}_r-Y^{\e,t,\xi}_{s+l}|^2I_{\{\sup\limits_{s\in[t,T]}|Y^{\e,t,\xi}_s|<R\}}\\
&\leq&CR\left(l^{1/2}+\(\mE(|K^{\e,t,\xi}|_t^{s+l}-|K^{\e,t,\xi}|_t^s)\)^{1/2}+\(\mE(|K^{\e,t,\xi}|_t^{s+l}-|K^{\e,t,\xi}|_t^s)^2\)^{1/2}\right).
\de
 
Besides, by the similar deduction to that for $(7)$ in \cite{pr1}, it holds that
 \ce
 \lim\limits_{l\rightarrow 0}\sup\limits_{\e}\mE(|K^{\e,t,\xi}|_t^{s+l}-|K^{\e,t,\xi}|_t^s)^2=0.
 \de

Finally, collecting the above estimates, we obtain the required limit by letting $l\rightarrow 0$ first and then $R\rightarrow\infty$. The proof is complete.
\end{proof}

Next, we estimate $\int_s^Tf(\frac{r}{\e},X_{r}^{\e,t,\xi},Y_{r}^{\e,t,\xi})\dif r-\int_s^T\bar{f}(X_{r}^{\e,t,\xi},Y_{r}^{\e,t,\xi})\dif r$.

\bl\label{fbarfcon}
Under the assumptions of Theorem \ref{averprin}, it holds that for $t\leq s\leq T$
\ce
\lim\limits_{\e\rightarrow0}\mE\left[\left|\int_s^Tf(\frac{r}{\e},X_{r}^{\e,t,\xi},Y_{r}^{\e,t,\xi})\dif r-\int_s^T\bar{f}(X_{r}^{\e,t,\xi},Y_{r}^{\e,t,\xi})\dif r\right|\right]=0.
\de
\el
\begin{proof}
First of all, we notice that
\ce
&&\mE\left|\int_s^Tf(\frac{r}{\e},X_{r}^{\e,t,\xi},Y_{r}^{\e,t,\xi})\dif r-\int_s^T\bar{f}(X_{r}^{\e,t,\xi},Y_{r}^{\e,t,\xi})\dif r\right|\\
&\leq&\mE\left|\int_s^T\left(f(\frac{r}{\e},X_{r}^{\e,t,\xi},Y_{r}^{\e,t,\xi})-f(\frac{r}{\e},X_{r(\d)}^{\e,t,\xi},Y_{r(\d)}^{\e,t,\xi})\right)\dif r\right|\\
&&+\mE\left|\int_s^T\left(f(\frac{r}{\e},X_{r(\d)}^{\e,t,\xi},Y_{r(\d)}^{\e,t,\xi})-\bar{f}(X_{r(\d)}^{\e,t,\xi},Y_{r(\d)}^{\e,t,\xi})\right)\dif r\right|\\
&&+\mE\left|\int_s^T\left(\bar{f}(X_{r(\d)}^{\e,t,\xi},Y_{r(\d)}^{\e,t,\xi})-\bar{f}(X_{r}^{\e,t,\xi},Y_{r}^{\e,t,\xi})\right)\dif r\right|\\
&=:& J_1+J_2+J_3.
\de

For $J_1+J_3$, by $(\mathbf{H}^1_{f,g})$, $(\mathbf{H}^{2}_{f})$ and Remark \ref{barfcond}, it holds that
\ce
J_1+J_3&\leq&L^{1/2}_3\mE\int_s^T(|X_{r}^{\e,t,\xi}-X_{r(\d)}^{\e,t,\xi}|+|Y_{r}^{\e,t,\xi}-Y_{r(\d)}^{\e,t,\xi}|)\dif r\\
&\leq&(2TL_3)^{1/2}\left(\int_s^T(\mE|X_{r}^{\e,t,\xi}-X_{r(\d)}^{\e,t,\xi}|^2+\mE|Y_{r}^{\e,t,\xi}-Y_{r(\d)}^{\e,t,\xi}|^2)\dif r\right)^{1/2}\\
&\leq&(2L_3)^{1/2}T\(\sup\limits_{v\in[t,T]}\mE\sup\limits_{v\leq r\leq v+\d}|X^{\e,t,\xi}_r-X^{\e,t,\xi}_v|^2\\
&&\qquad\qquad\qquad +\sup\limits_{v\in[t,T]}\sup\limits_{v\leq r\leq v+\d}\mE|Y^{\e,t,\xi}_r-Y^{\e,t,\xi}_v|^2\)^{1/2}\\
&\leq&(2L_3)^{1/2}T\(\sup\limits_{\e}\sup\limits_{v\in[t,T]}\mE\sup\limits_{v\leq r\leq v+\d}|X^{\e,t,\xi}_r-X^{\e,t,\xi}_v|^2\\
&&\qquad\qquad\qquad+\sup\limits_{\e}\sup\limits_{v\in[t,T]}\sup\limits_{v\leq r\leq v+\d}\mE|Y^{\e,t,\xi}_r-Y^{\e,t,\xi}_v|^2\)^{1/2},
\de
where $r(\d):=[\frac{r-s}{\d}]\d+s$ for any $0<\d<T-s$.

For $J_2$, by the similar deduction to that for $I_2$ in the proof of Lemma \ref{xeconver}, one can have that
\ce
&&\int_s^T\left(f(\frac{r}{\e},X_{r(\d)}^{\e,t,\xi},Y_{r(\d)}^{\e,t,\xi})-\bar{f}(X_{r(\d)}^{\e,t,\xi},Y_{r(\d)}^{\e,t,\xi})\right)\dif r\\
&=&\sum\limits_{k=0}^{[\frac{T-s}{\d}]-1}\e\int_{\frac{k\d+s}{\e}}^{\frac{(k+1)\d+s}{\e}}\left(f(\tau, X^{\e,t,\xi}_{k\d+s},Y^{\e,t,\xi}_{k\d+s})-\bar{f}(X^{\e,t,\xi}_{k\d+s},Y^{\e,t,\xi}_{k\d+s})\right)\dif \tau\\
&&+\int_{[\frac{T-s}{\d}]\d}^{T-s}\left(f(\frac{\varrho+s}{\e}, X^{\e,t,\xi}_{[\frac{\varrho}{\d}]\d+s},Y^{\e,t,\xi}_{[\frac{\varrho}{\d}]\d+s})-\bar{f}(X^{\e,t,\xi}_{[\frac{\varrho}{\d}]\d+s},Y^{\e,t,\xi}_{[\frac{\varrho}{\d}]\d+s})\right)\dif \varrho.
\de
So, $(\mathbf{H}^1_{f,g})$, $(\mathbf{H}^{2}_{f})$, (\ref{xeboun}), (\ref{yees}) and Remark \ref{barfcond} imply that
\ce
J_2&\leq&\mE\sum\limits_{k=0}^{[\frac{T-s}{\d}]-1}\Bigg|\e\int_{\frac{k\d+s}{\e}}^{\frac{(k+1)\d+s}{\e}}\left(f(\tau, X^{\e,t,\xi}_{k\d+s},Y^{\e,t,\xi}_{k\d+s})-\bar{f}(X^{\e,t,\xi}_{k\d+s},Y^{\e,t,\xi}_{k\d+s})\right)\dif \tau\Bigg|\\
&&+\mE\int_{[\frac{T-s}{\d}]\d}^{T-s}\left|f(\frac{\varrho+s}{\e}, X^{\e,t,\xi}_{[\frac{\varrho}{\d}]\d+s},Y^{\e,t,\xi}_{[\frac{\varrho}{\d}]\d+s})-\bar{f}(X^{\e,t,\xi}_{[\frac{\varrho}{\d}]\d+s},Y^{\e,t,\xi}_{[\frac{\varrho}{\d}]\d+s})\right|\dif \varrho\\
&\leq&\mE\sum\limits_{k=0}^{[\frac{T-s}{\d}]-1}\Bigg|\e\int_{0}^{\frac{(k+1)\d+s}{\e}}\left(f(\tau, X^{\e,t,\xi}_{k\d+s},Y^{\e,t,\xi}_{k\d+s})-\bar{f}(X^{\e,t,\xi}_{k\d+s},Y^{\e,t,\xi}_{k\d+s})\right)\dif \tau\\
&&\qquad\qquad\quad-\e\int_{0}^{\frac{k\d+s}{\e}}\left(f(\tau, X^{\e,t,\xi}_{k\d+s},Y^{\e,t,\xi}_{k\d+s})-\bar{f}(X^{\e,t,\xi}_{k\d+s},Y^{\e,t,\xi}_{k\d+s})\right)\dif \tau\Bigg|\\
&&+C\d\left(\mE(1+\sup\limits_{r\in[t,T]}|X^{\e,t,\xi}_r|^2+\sup\limits_{r\in[t,T]}|Y^{\e,t,\xi}_r|^2)\right)^{1/2}\\
&\leq&\mE\sum\limits_{k=0}^{[\frac{T-s}{\d}]-1}\Bigg|\e\int_{0}^{\frac{(k+1)\d+s}{\e}}\left(f(\tau, X^{\e,t,\xi}_{k\d+s},Y^{\e,t,\xi}_{k\d+s})-\bar{f}(X^{\e,t,\xi}_{k\d+s},Y^{\e,t,\xi}_{k\d+s})\right)\dif \tau\\
&&\qquad\qquad\quad-\e\int_{0}^{\frac{k\d+s}{\e}}\left(f(\tau, X^{\e,t,\xi}_{k\d+s},Y^{\e,t,\xi}_{k\d+s})-\bar{f}(X^{\e,t,\xi}_{k\d+s},Y^{\e,t,\xi}_{k\d+s})\right)\dif \tau\Bigg|\\
&&+C\d\left(1+\mE|\xi|^2+C\right)^{1/2}.
\de

Combining the above deduction, we obtain that
\ce
0&\leq& \mE\left|\int_s^Tf(\frac{r}{\e},X_{r}^{\e,t,\xi},Y_{r}^{\e,t,\xi})\dif r-\int_s^T\bar{f}(X_{r}^{\e,t,\xi},Y_{r}^{\e,t,\xi})\dif r\right|\\
&\leq&(2L_3)^{1/2}T\(\sup\limits_{\e}\sup\limits_{v\in[t,T]}\mE\sup\limits_{v\leq r\leq v+\d}|X^{\e,t,\xi}_r-X^{\e,t,\xi}_v|^2\\
&&\qquad\qquad\qquad+\sup\limits_{\e}\sup\limits_{v\in[t,T]}\sup\limits_{v\leq r\leq v+\d}\mE|Y^{\e,t,\xi}_r-Y^{\e,t,\xi}_v|^2\)^{1/2}\\
&&+\mE\sum\limits_{k=0}^{[\frac{T-s}{\d}]-1}\Bigg|\e\int_{0}^{\frac{(k+1)\d+s}{\e}}\left(f(\tau, X^{\e,t,\xi}_{k\d+s},Y^{\e,t,\xi}_{k\d+s})-\bar{f}(X^{\e,t,\xi}_{k\d+s},Y^{\e,t,\xi}_{k\d+s})\right)\dif \tau\\
&&\qquad\qquad\quad-\e\int_{0}^{\frac{k\d+s}{\e}}\left(f(\tau, X^{\e,t,\xi}_{k\d+s},Y^{\e,t,\xi}_{k\d+s})-\bar{f}(X^{\e,t,\xi}_{k\d+s},Y^{\e,t,\xi}_{k\d+s})\right)\dif \tau\Bigg|\\
&&+C\d\left(1+\mE|\xi|^2+C\right)^{1/2}.
\de
Letting $\e\rightarrow 0$ first and then $\d\rightarrow 0$, by $(\mathbf{H}^{2}_{f})$, the dominated convergence theorem and Lemma \ref{yessles}, we get the required convergence.
\end{proof}

Now we are going to prove Theorem \ref{averprin}.

{\bf Proof of Theorem \ref{averprin}.} First of all, by Lemma \ref{8proctigh}, Theorem \ref{xbarxdiff} and Lemma \ref{stopotigh}, we know that there exists a subsequence $\left\{\e_k\right\}_{k\in\mN}$ decreasing to $0$, four processes $\hat{Y}^{t,\xi}, \hat{M}^{t,\xi}, \hat{{\bf U}}, \hat{{\bf V}}$ in $D([t,T],\mR^d)$, and a countable set $I \subset[t, T)$ such that $(X^{\e_k,t,\xi}, K^{\e_k,t,\xi}, Y^{\e_k,t,\xi}, \\M^{\e_k,t,\xi}, {\bf U}^{\e_k}, {\bf V}^{\e_k})$ converges in distribution to $(\bar{X}^{t,\xi}, \bar{K}^{t,\xi}, \hat{Y}^{t,\xi}, \hat{M}^{t,\xi}, \hat{{\bf U}}, \hat{{\bf V}})$ on $[t, T] \backslash I$. 

We observe the limits of $\int_s^Tg(r,X_{r}^{\e_k,t,\xi},Y_{r}^{\e_k,t,\xi})\dif |K^{\e_k,t,\xi}|_t^{r}$ and $\int_s^Tf(\frac{r}{\e_k},X_{r}^{\e_k,t,\xi},Y_{r}^{\e_k,t,\xi})\dif r$ as $k\rightarrow\infty$. On one hand, $(\mathbf{H}^1_{f,g})$ and Lemma \ref{finivarilimi} imply that for any $t\leq s\leq T$, as $k\rightarrow\infty$
\ce
\int_s^Tg(r,X_{r}^{\e_k,t,\xi},Y_{r}^{\e_k,t,\xi})\dif |K^{\e_k,t,\xi}|_t^{r}&=&\int_s^Tg(r,X_{r}^{\e_k,t,\xi},Y_{r}^{\e_k,t,\xi})\<\triangledown \phi(X_{r}^{\e_k,t,\xi}), \dif K_{r}^{\e_k,t,\xi}\>\\
&\overset{in~distribution}{\longrightarrow}& \int_s^Tg(r,\bar{X}_r^{t,\xi},\hat{Y}_r^{t,\xi})\<\triangledown \phi(\bar{X}_{r}^{t,\xi}), \dif \bar{K}_{r}^{t,\xi}\>\\
&=&\int_s^Tg(r,\bar{X}_r^{t,\xi},\hat{Y}_r^{t,\xi})\dif |\bar{K}^{t,\xi}|_t^{r}.
\de
On the other hand, by Lemma \ref{fbarfcon}, (\ref{xeboun}) and (\ref{yees}), it holds that as $k\rightarrow\infty$
\ce
\int_s^Tf(\frac{r}{\e_k},X_{r}^{\e_k,t,\xi},Y_{r}^{\e_k,t,\xi})\dif r\overset{in~distribution}{\longrightarrow}\int_s^T\bar{f}(\bar{X}_{r}^{t,\xi},\hat{Y}_{r}^{t,\xi})\dif r.
\de

Now we take the limit in (\ref{yesatiequa}) and obtain that for any $s\in [t, T] \backslash I$
\be
\hat{Y}^{t,\xi}_s&=&\Phi(\bar{X}_{T}^{t,\xi})-(\hat{{\bf U}}_T-\hat{{\bf U}}_s)-(\hat{{\bf V}}_T-\hat{{\bf V}}_s)+\int_s^T\bar{f}(\bar{X}_{r}^{t,\xi},\hat{Y}_{r}^{t,\xi})\dif r\no\\
&&+\int_s^Tg(r,\bar{X}_r^{t,\xi},\hat{Y}_r^{t,\xi})\dif |\bar{K}^{t,\xi}|_t^{r}-(\hat{M}^{t,\xi}_T-\hat{M}^{t,\xi}_s).
\label{hatyequa}
\ee
Since $\hat{Y}^{t,\xi}, \hat{M}^{t,\xi}, \hat{{\bf U}}, \hat{{\bf V}}$ are c\`adl\`ag in $s$, the above equation holds for any $s\in [t,T]$. Moreover, from the above equation, it follows that 
$\hat{M}^{t,\xi}$ is $(\sF^{B,\hat{Y}^{t,\xi}, \hat{M}^{t,\xi}, \hat{{\bf U}}, \hat{{\bf V}}}_s)_{s\in[t,T]}$-adapted.

Next,  we know that $M^{\e_k,t,\xi}$ is a martingales with respect to $(\sF^{B,Y^{\e_k,t,\xi},M^{\e_k,t,\xi}, {\bf U}^{\e_k}, {\bf V}^{\e_k}}_s)_{s\in[t,T]}$. Furthermore, by (\ref{zees}), it holds that
\ce
\sup\limits_{\e}\mE[\sup\limits_{s\in[t,T]}|M^{\e_k,t,\xi}_{s}|^2]\leq 4\sup\limits_{\e}\mE\int_t^T\|Z_{r}^{\e,t,\xi}\|^2\dif r\leq C.
\de
Thus, by Lemma \ref{martsto}, we have that $\hat{M}^{t,\xi}$ is a martingales with respect to $(\sF^{B,\hat{Y}^{t,\xi}, \hat{M}^{t,\xi}, \hat{{\bf U}}, \hat{{\bf V}}}_s)_{s\in[t,T]}$.

In the following, we notice that $(\bar{Y}^{t,\xi},\bar{Z}^{t,\xi},\bar{U}^{t,\xi}, \bar{V}^{t,\xi})$ solves uniquely the following equation
\be
\bar{Y}^{t,\xi}_s&=&\Phi(\bar{X}_{T}^{t,\xi})-(\bar {\bf U}_T-\bar {\bf U}_s)-(\bar {\bf V}_T-\bar {\bf V}_s)+\int_s^T\bar{f}(\bar{X}_{r}^{t,\xi},\bar{Y}_{r}^{t,\xi})\dif r\no\\
&&+\int_s^Tg(r,\bar X_{r}^{t,\xi},\bar Y_{r}^{t,\xi})\dif |\bar K^{t,\xi}|_t^r-(\bar M^{t,\xi}_T-\bar M^{t,\xi}_s),
\label{baryequa}
\ee
where 
$$
\bar {\bf U}_s:=\int_t^s\bar{U}^{t,\xi}_r\dif r, \quad \bar {\bf V}_s:=\int_t^s\bar{V}^{t,\xi}_r\dif |\bar{K}^{t,\xi}|_t^r, \quad \bar M^{t,\xi}_s:=\int_t^s \bar Z_{r}^{t,\xi}\dif B_{r}.
$$
Since $\mE\int_t^T \|\bar Z_{r}^{t,\xi}\|^2\dif r<\infty$ and $B$ is $(\sF^{B,\hat{Y}^{t,\xi}, \hat{M}^{t,\xi}, \hat{{\bf U}}, \hat{{\bf V}}}_s)_{s\in[t,T]}$-Brownian motion, $\bar M^{t,\xi}$ is an $(\sF^{B,\hat{Y}^{t,\xi}, \hat{M}^{t,\xi}, \hat{{\bf U}}, \hat{{\bf V}}}_s)_{s\in[t,T]}$-martingale.

By the It\^o formula, it holds that for any two stopping times $t\leq \tau_1<\tau_2\leq T$ $\mP.$-a.s.
\be
&&|\hat{Y}^{t,\xi}_{\tau_1}-\bar{Y}^{t,\xi}_{\tau_1}|^2+\int_{\tau_1}^{\tau_2}\dif [\hat{M}^{t,\xi}+\hat{{\bf U}}+\hat{{\bf V}}-\bar M^{t,\xi}-\bar {\bf U}-\bar {\bf V}]_r\no\\
&&+2\int_{\tau_1}^{\tau_2}\<\hat{Y}^{t,\xi}_r-\bar{Y}^{t,\xi}_r, \dif (\hat{{\bf U}}_r+\hat{{\bf V}}_r)-\dif (\bar {\bf U}_r+\bar {\bf V}_r)\>\no\\
&=&|\hat{Y}^{t,\xi}_{\tau_2}-\bar{Y}^{t,\xi}_{\tau_2}|^2+2\int_{\tau_1}^{\tau_2}\<\hat{Y}^{t,\xi}_r-\bar{Y}^{t,\xi}_r, \bar{f}(\bar{X}_{r}^{t,\xi},\hat{Y}_{r}^{t,\xi})-\bar{f}(\bar{X}_{r}^{t,\xi},\bar{Y}_{r}^{t,\xi})\>\dif r\no\\
&&+2\int_{\tau_1}^{\tau_2}\<\hat{Y}^{t,\xi}_r-\bar{Y}^{t,\xi}_r, g(r,\bar{X}_r^{t,\xi},\hat{Y}_r^{t,\xi})-g(r,\bar X_{r}^{t,\xi},\bar Y_{r}^{t,\xi})\>\dif |\bar{K}^{t,\xi}|_t^{r}\no\\
&&-2\int_{\tau_1}^{\tau_2}\<\hat{Y}^{t,\xi}_r-\bar{Y}^{t,\xi}_r,\dif (\hat{M}^{t,\xi}_r-\bar M^{t,\xi}_r)\>,
\label{itohatybary}
\ee
where $[\hat{M}^{t,\xi}+\hat{{\bf U}}+\hat{{\bf V}}-\bar M^{t,\xi}-\bar {\bf U}-\bar {\bf V}]$ is the quadratic variation process of $\hat{M}^{t,\xi}+\hat{{\bf U}}+\hat{{\bf V}}-\bar M^{t,\xi}-\bar {\bf U}-\bar {\bf V}$. Set
\ce
&&\bar{X}_{s}^{t,\xi}=\xi, |\bar{K}^{t,\xi}|_t^s=0, \hat{Y}_{s}^{t,\xi}=\hat{Y}_{t}^{t,\xi}, \hat{M}_{s}^{t,\xi}=0, \hat{{\bf U}}_s=0, \hat{{\bf V}}_s=0, \quad 0\leq s<t,\\
&&\bar{X}_{s}^{t,\xi}=\bar{X}_{T}^{t,\xi}, |\bar{K}^{t,\xi}|_t^s=|\bar{K}^{t,\xi}|_t^T, \hat{Y}_{s}^{t,\xi}=\hat{Y}_{T}^{t,\xi}, \hat{M}_{s}^{t,\xi}=\hat{M}_{T}^{t,\xi}, \hat{{\bf U}}_s=\hat{{\bf U}}_T, \hat{{\bf V}}_s=\hat{{\bf V}}_T, \quad s\geq T,
\de
and
\ce
&&\bar{Y}_{s}^{t,\xi}=\bar{Y}_{t}^{t,\xi}, \bar{M}_{s}^{t,\xi}=0, \bar{{\bf U}}_s=0, \bar{{\bf V}}_s=0, \quad 0\leq s<t,\\
&&\bar{Y}_{s}^{t,\xi}=\bar{Y}_{T}^{t,\xi}, \bar{M}_{s}^{t,\xi}=\bar{M}_{T}^{t,\xi}, \bar{{\bf U}}_s=\bar{{\bf U}}_T, \bar{{\bf V}}_s=\bar{{\bf V}}_T, \quad s\geq T,
\de
and we extend (\ref{itohatybary}) to any two stopping times $0\leq \tau_1<\tau_2<\infty$ $\mP.$-a.s. Thus, from \cite[Lemma 7]{mr2} and $(\mathbf{H}^1_{f,g})$, $(\mathbf{H}^{2}_{f})$, it follows that
\be
&&\mE|\hat{Y}^{t,\xi}_{\tau_1}-\bar{Y}^{t,\xi}_{\tau_1}|^2+\mE\left([\hat{M}^{t,\xi}-\bar M^{t,\xi}]_{\tau_2}-[\hat{M}^{t,\xi}-\bar M^{t,\xi}]_{\tau_1}\right)\no\\
&\leq&\mE|\hat{Y}^{t,\xi}_{\tau_2}-\bar{Y}^{t,\xi}_{\tau_2}|^2+2(L_3^{1/2}+L_4^{1/2})\mE\int_{\tau_1}^{\tau_2}|\hat{Y}^{t,\xi}_r-\bar{Y}^{t,\xi}_r|^2\dif (r+|\bar{K}^{t,\xi}|_t^{r}).
\label{hatybarymm}
\ee

Set $\kappa_s(\omega):=s+|\bar{K}^{t,\xi}|_t^{s \wedge T}$, and $s\mapsto\kappa_s$ is a continuous strictly increasing and bijective function and $\kappa^{-1}$ denotes its inverse mapping. Let us consider the stopping times $\tau_1=\kappa_{s_1}^{-1}$ and $\tau_2=\kappa_{s_2}^{-1}$, where $0 \leq s_1 < s_2$ are non random. From (\ref{hatybarymm}), it follows that
$$
\mathbb{E}\left|\hat{Y}^{t,\xi}_{\kappa_{s_1}^{-1}}-\bar{Y}^{t,\xi}_{\kappa_{s_1}^{-1}}\right|^2 \leq \mathbb{E}\left|\hat{Y}^{t,\xi}_{\kappa_{s_2}^{-1}}-\bar{Y}^{t,\xi}_{\kappa_{s_2}^{-1}}\right|^2+2(L_3^{1/2}+L_4^{1/2})\int_{s_1}^{s_2} \mathbb{E}\left|\hat{Y}^{t,\xi}_{\kappa_{r}^{-1}}-\bar{Y}^{t,\xi}_{\kappa_{r}^{-1}}\right|^2 \dif r,
$$
which together with the Gronwall inequality yields that
$$
\mathbb{E}e^{2(L_3^{1/2}+L_4^{1/2}) s_1}\left|\hat{Y}^{t,\xi}_{\kappa_{s_1}^{-1}}-\bar{Y}^{t,\xi}_{\kappa_{s_1}^{-1}}\right|^2 \leq \mathbb{E}e^{2(L_3^{1/2}+L_4^{1/2}) s_2}\left|\hat{Y}^{t,\xi}_{\kappa_{s_2}^{-1}}-\bar{Y}^{t,\xi}_{\kappa_{s_2}^{-1}}\right|^2.
$$
Note that
$$
e^{2(L_3^{1/2}+L_4^{1/2}) s}\left|\hat{Y}^{t,\xi}_{\kappa_{s}^{-1}}-\bar{Y}^{t,\xi}_{\kappa_{s}^{-1}}\right|^2=0, \quad \text { for any } s \geq \kappa_T, \text { a.s. }
$$
and
$$
\sup _{s \geq 0}e^{2(L_3^{1/2}+L_4^{1/2}) s}\left|\hat{Y}^{t,\xi}_{\kappa_{s}^{-1}}-\bar{Y}^{t,\xi}_{\kappa_{s}^{-1}}\right|^2 \leq \sup _{r \in[t, T]}e^{2(L_3^{1/2}+L_4^{1/2}) \kappa_r}\left|\hat{Y}^{t,\xi}_r-\bar{Y}^{t,\xi}_r\right|^2.
$$
Besides, by the similar deduction to that in \cite[Theorem 9]{mr1}, it holds that for any $\mu\geq 0$, there exists a constant $C$ such that
\ce
\mE \sup\limits_{s\in[t,T]}e^{\mu |\bar{K}^{t,\xi}|_t^s}(|\hat{Y}^{t,\xi}_s|^2+|\bar{Y}^{t,\xi}_s|^2)\leq C.
\de
As $s_2 \rightarrow \infty$, by the dominated convergence theorem we deduce that
$$
\mathbb{E}e^{2(L_3^{1/2}+L_4^{1/2}) s_1}\left|\hat{Y}^{t,\xi}_{\kappa_{s_1}^{-1}}-\bar{Y}^{t,\xi}_{\kappa_{s_1}^{-1}}\right|^2=0
$$
for any $s_1 \geq 0$, which yields that
\ce
\hat{Y}^{t,\xi}_s=\bar{Y}^{t,\xi}_s, \quad s\in[t,T], \quad \mP-a.s..
\de
and by (\ref{hatybarymm})
\ce
\hat{M}^{t,\xi}_s=\bar M^{t,\xi}_s, \quad s\in[t,T], \quad \mP-a.s..
\de

Finally, combining the above deduction with (\ref{hatyequa}) and (\ref{baryequa}), we conclude that 
\ce
\hat{{\bf U}}_s+\hat{{\bf V}}_s=\bar{{\bf U}}_s+\bar{{\bf V}}_s, \quad s\in[t,T], \quad \mP-a.s..
\de
The proof is complete.

{\bf Proof of Corollary \ref{xaverprin}.} Since $Y_{t}^{\e,t,x}$ and $\bar{Y}^{t,x}_t$ are deterministic, we know that
\ce
Y_{t}^{\e,t,x}&=&\mE\bigg[\Phi(X_{T}^{\e,t,x})-({\bf U}^{\e}_T+{\bf V}^{\e}_T)+\int_t^Tf(\frac{r}{\e},X_{r}^{\e,t,x},Y_{r}^{\e,t,x})\dif r\\
&&\qquad+\int_t^Tg(r,X_{r}^{\e,t,x},Y_{r}^{\e,t,x})\dif |K^{\e,t,x}|_t^r\bigg]
\de
and
\ce
\bar{Y}^{t,x}_t=\mE\left[\Phi(\bar X_{T}^{t,x})-(\bar{\bf U}_T+\bar{\bf V}_T)+\int_t^T\bar{f}(\bar{X}_{r}^{t,x},\bar{Y}_{r}^{t,x})\dif r+\int_t^Tg(r,\bar X_{r}^{t,x},\bar Y_{r}^{t,x})\dif |\bar K^{t,x}|_t^r\right].
\de

Next, we set 
\ce
&&\cK^\e:=\Phi(X_{T}^{\e,t,x})-({\bf U}^{\e}_T+{\bf V}^{\e}_T)+\int_t^Tf(\frac{r}{\e},X_{r}^{\e,t,x},Y_{r}^{\e,t,x})\dif r+\int_t^Tg(r,X_{r}^{\e,t,x},Y_{r}^{\e,t,x})\dif |K^{\e,t,x}|_t^r,\\
&&\bar\cK:=\Phi(\bar X_{T}^{t,x})-(\bar{\bf U}_T+\bar{\bf V}_T)+\int_t^T\bar{f}(\bar{X}_{r}^{t,x},\bar{Y}_{r}^{t,x})\dif r+\int_t^Tg(r,\bar X_{r}^{t,x},\bar Y_{r}^{t,x})\dif |\bar K^{t,x}|_t^r,
\de
and the H\"older inequality and $(\mathbf{H}^1_{f,g})$ yield that
\ce
\mE|\cK^\e|^2&\leq& 5\mE|\Phi(X_{T}^{\e,t,x})|^2+5\mE|{\bf U}^{\e}_T|^2+5\mE|{\bf V}^{\e}_T|^2+5\mE\left|\int_t^Tf(\frac{r}{\e},X_{r}^{\e,t,x},Y_{r}^{\e,t,x})\dif r\right|^2\\
&&+5\mE\left|\int_t^Tg(r,X_{r}^{\e,t,x},Y_{r}^{\e,t,x})\dif |K^{\e,t,x}|_t^r\right|^2\\
&\leq&5C+5T\mE\int_t^T|U^{\e,t,\xi}_r|^2\dif r+5\left(\mE\int_t^T|V^{\e,t,\xi}_r|^2\dif |K^{\e,t,\xi}|_t^r\right)\left(\mE|K^{\e,t,\xi}|_t^T\right)\\
&&+C\left(\mE\int_t^T(1+|Y_{r}^{\e,t,\xi}|)^2(\dif r+\dif |K^{\e,t,\xi}|_t^r)\right)\left(T+\mE|K^{\e,t,\xi}|_t^T\right)\\
&&+C\left(\mE\sup\limits_{r\in[t,T]}|X^{\e,t,\xi}_r|^2\right)\left(T^2+\mE(|K^{\e,t,\xi}|_t^T)^2\right).
\de
Thus, (\ref{uevees}), (\ref{yees}), (\ref{keboun}), (\ref{xeboun}) and \cite[Proposition 1]{mr2} assure that
\ce
\sup\limits_{\e}\mE|\cK^\e|^2\leq C,
\de
which implies that $\{\cK^\e\}$ is uniformly integrable. Besides, by the proof of Theorem \ref{averprin}, it holds that $\cK^\e$ converges in distribution to $\bar \cK$.

Finally, combining the above deduction, we have that $\mE\cK^\e$ converges to $\mE\bar\cK$. That is, $Y_{t}^{\e,t,x}$ converges to $\bar{Y}^{t,x}_t$. The proof is complete.

\section{Application}\label{app}

In this section, we apply our result to a type of multivalued Dirichlet-Neumann problems and study its homogenization.

First of all, we consider the system (\ref{dnpdein}) with the functions $\varphi, \psi: \mR^d\rightarrow (-\infty, \infty]$ decoupled in the sense that $\varphi\left(u_1, \ldots, u_d\right)=\varphi_1\left(u_1\right)+\cdots+\varphi_d\left(u_d\right)$ and $\psi\left(u_1, \ldots, u_d\right)=\psi_1\left(u_1\right)+$ $\cdots+\psi_d\left(u_d\right)$, where $\left.\left.\varphi_i, \psi_i: \mathbb{R} \rightarrow\right (-\infty,+\infty\right]$ are lower semicontinuous convex functions; hence $\partial \varphi\left(u_1, \ldots, u_d\right)=\partial \varphi_1\left(u_1\right) \times \cdots \times \partial \varphi_d\left(u_d\right)$ and similarly for $\partial \psi$.

Consider two following systems of parabolic variation inequalities
\be\left\{\begin{array}{l}
\frac{\p u_i^\e(t,x)}{\p t}+\sL^\e u_i^\e(t,x)+f_i(\frac{t}{\e},x,u^\e(t,x))\in \p \varphi_i(u_i^\e(t,x)), ~(t,x)\in[0,T]\times\cO, i=1, \cdots, d\\ 
\frac{\p u_i^\e(t,x)}{\p n}+g_i(t,x,u^\e(t,x))\in \p\psi_i(u_i^\e(t,x)), ~(t,x)\in[0,T]\times\p \cO, i=1, \cdots, d\\
u^\e(T,x)=\Phi(x), \quad x\in\bar{\cO},
\end{array}
\right.
\label{dpdein1}
\ee
and
\be\left\{\begin{array}{l}
\frac{\p u_i(t,x)}{\p t}+(\bar \sL u_i)(t,x)+\bar{f}_i(x,u(t,x))\in \p \varphi_i(u_i(t,x)), ~(t,x)\in[0,T]\times\cO, i=1, \cdots, d\\ 
\frac{\p u_i(t,x)}{\p n}+g_i(t,x,u(t,x))\in \p\psi_i(u_i(t,x)), ~(t,x)\in[0,T]\times\p \cO, i=1, \cdots, d\\
u(T,x)=\Phi(x), \quad x\in\bar{\cO},
\end{array}
\right.
\label{dpdein2}
\ee
where $\sL^\e$ and $\bar\sL$ are defined in (\ref{ledefi}) and (\ref{barldefi}), respectively. The following theorem describes the relationship between the system (\ref{dpdein1}) and the system (\ref{dpdein2}).

\bt\label{dnbounth}
Assume that $(\mathbf{H}^{1}_{b, \s})$, $(\mathbf{H}^{2}_{b, \s})$, $(\mathbf{H}^3_{\s})$, $({\bf H}_{\varphi,\psi})$, $(\mathbf{H}^1_{f,g})$ and $(\mathbf{H}^{2}_{f})$ hold. Then $u^\e(t,x)$ converges to $u(t,x)$ as $\e$ tends to $0$, where $u^\e(t,x)$ and $u(t,x)$ are  the unique viscosity solutions  of the system (\ref{dpdein1}) and the system (\ref{dpdein2}), respectively.
\et
\begin{proof}
By the nonlinear Feynman-Kac formula (\cite[Theorem 5.2]{pr1}), we know that $u^\e(t,x):=Y_{t}^{\e,t,x}$ and $u(t,x):=\bar{Y}_t^{t,x}$  are unique viscosity solutions of the system (\ref{dpdein1}) and the system (\ref{dpdein2}), respectively, where $Y^{\e,t,x}$ solves the following forward-backward coupled stochastic variational inequalities
\ce\left\{\begin{array}{l}
\dif X_{s}^{\e,t,x}=\triangledown \phi(X_{s}^{\e,t,x})\dif |K^{\e,t,x}|_t^s+b(\frac{s}{\e},X_{s}^{\e,t,x})\dif s+\s(\frac{s}{\e},X_{s}^{\e,t,x})\dif B_{s},\\
K^{\e,t,x}_s=\int_t^s\triangledown \phi(X_{r}^{\e,t,x})\dif |K^{\e,t,x}|_t^r, \quad |K^{\e,t,x}|_t^s=\int_t^s I_{\{X_{r}^{\e,t,x}\in \p\cO\}}\dif |K^{\e,t,x}|_t^r,\\
X_{t}^{\e,t,x}=x\in \overline{\cO},\\
\dif Y_{s}^{\e,t,x}\in\p\varphi(Y_{s}^{\e,t,x})\dif s+\p\psi(Y_{s}^{\e,t,x})\dif |K^{\e,t,x}|_t^s-f(\frac{s}{\e},X_{s}^{\e,t,x},Y_{s}^{\e,t,x})\dif s\\
\qquad\qquad-g(s,X_{s}^{\e,t,x},Y_{s}^{\e,t,x})\dif |K^{\e,t,x}|_t^s+Z_{s}^{\e,t,x}\dif B_{s},\\
Y_{T}^{\e,t,x}=\Phi(X_{T}^{\e,t,x}).
\end{array}
\right.
\de
and $\bar{Y}^{t,x}$ solves the following forward-backward coupled stochastic variational inequalities
\ce\left\{\begin{array}{l}
\dif \bar{X}_{s}^{t,x}=\triangledown \phi(\bar{X}_{s}^{t,x})\dif |\bar K^{t,x}|_t^s+\bar{b}(\bar{X}_{s}^{t,x})\dif s+\bar{\s}(\bar{X}_{s}^{t,x})\dif B_{s},\\
\bar K^{t,x}_s=\int_t^s\triangledown \phi(\bar X_{r}^{,x})\dif |\bar K^{t,x}|_t^r, \quad |\bar K^{t,x}|_t^s=\int_t^s I_{\{\bar X_{r}^{t,x}\in \p\cO\}}\dif |\bar K^{t,x}|_t^r,\\
\bar{X}_{t}^{t,x}=x\in \overline{\cO},\\
\dif \bar{Y}_{s}^{t,x}\in \p\varphi(\bar{Y}_{s}^{t,x})\dif s+\p\psi(\bar{Y}_{s}^{t,x})\dif |\bar{K}^{t,x}|_t^s-\bar{f}(\bar{X}_{s}^{t,x},\bar{Y}_{s}^{t,x})\dif s\\
\qquad\qquad-g(s,\bar{X}_{s}^{t,x},\bar{Y}_{s}^{t,x})\dif |\bar{K}^{t,x}|_t^s+\bar{Z}_{s}^{t,x}\dif B_{s},\\
\bar{Y}_{T}^{t,x}=\Phi(\bar{X}_{T}^{t,x}).
\end{array}
\right.
\de
Moreover, by Corollary \ref{xaverprin}, it holds that $Y^{\e,t,x}_t$ converges to $\bar{Y}^{t,x}_t$ as $\e$ tends to $0$. Based on this, we conclude that $u^\e(t,x)$ converges to $u(t,x)$. That is, we establish homogenization for the system (\ref{dpdein1}). 
\end{proof}

\br
$(i)$ If $\cO=\mR^m$, that is $K^{\e,t,x}=0$, and $f(t,x,y)$ is independent of $t$, Theorem \ref{dnbounth} becomes \cite[Theorem 4.1]{eo}.

$(ii)$ If $\cO=\mR^m$ and $\varphi=\psi=0$, Theorem \ref{dnbounth} is just right \cite[Theorem 6.4]{q1} with $f_2(z)=0$.
\er

\end{document}